\documentclass{elsart}
\journal{Journal of Algebra}
\usepackage{amsfonts}
\usepackage{amscd}
\usepackage{amssymb}
\usepackage{amsmath}
\usepackage{latexsym}
\usepackage{diagrams}

\newtheorem{theorem}{Theorem}[section]
\newtheorem{lemma}[theorem]{Lemma}
\newtheorem{proposition}[theorem]{Proposition}
\newtheorem{corollary}[theorem]{Corollary}

\newtheorem{definition}[theorem]{Definition}

\newtheorem{remark}[theorem]{Remark}

\def\CQFD{$\Box$}
\newenvironment{proof}{{\it Proof}.}{\hfill \CQFD \\}

\def\PP{{\mathbb{P}}}
\def\AA{{\mathbb{A}}}

\def\NN{{\mathbb{N}}}
\def\ZZ{{\mathbb{Z}}}
\def\GG{{\mathbb{G}}}

\def\OO{{\mathcal{O}}}
\def\CC{{\mathcal{C}}}
\def\BB{{\mathcal{B}}}

\def\Zc{{\mathcal{Z}}}
\def\Bc{{\mathcal{B}}}
\def\Mc{{\mathcal{M}}}
\def\Lc{{\mathcal{L}}}

\def\fg{\mathbf{f}}
\def\Tg{\mathbf{T}}

\def\ag{\mathbf{a}}

\def\gr{\mathrm{gr}}
\def\dim{\mathrm{dim}}
\def\det{\mathrm{det}}
\def\div{\mathrm{div}}
\def\deg{\mathrm{deg}}
\def\Proj{\mathrm{Proj}}
\def\Biproj{\mathrm{Biproj}}

\def\pd{\mathrm{pd}}
\def\Spec{\mathrm{Spec}}
\def\codim{\mathrm{codim}}
\def\depth{\mathrm{depth}}
\def\Sym{\mathrm{Sym}}
\def\ker{\mathrm{ker}}
\def\ann{\mathrm{ann}}
\def\TF{\mathrm{TF}}
\def\Rees{\mathrm{Rees}}

\def\Im{\mathrm{Im}}
\def\H{\mathrm{H}}
\def\P{\mathrm{P}}
\def\mult{\mathrm{mult}}
\def\length{\mathrm{length}}
\def\Bez{{\mathbb{B}ez}}
\def\Res{{\mathrm{Res}}}

\def\mm{\mathfrak{m}}
\def\pp{\mathfrak{p}}

\begin{document}

\begin{frontmatter}
\title{On the closed image of a rational map and the implicitization problem}
\author{Laurent Bus\'e}
\address{Universit\'e de Nice Sophia-Antipolis, UMR 6621, 
   Parc Valrose, BP 71,
  06108, Nice Cedex 02, France}
\author{Jean-Pierre Jouanolou}
\address{Universit\'e Louis Pasteur, D\'epartement de Math\'ematiques, 
  7 rue Ren\'e Descartes,
  67084, Strasbourg Cedex, France}

\begin{abstract}
In this paper, we investigate some topics around the closed image $S$ 
of a rational map $\lambda$ given by some homogeneous elements $f_1,\ldots,f_n$
of the same degree in a graded algebra $A$. We first compute the degree 
of this closed image in case $\lambda$ is generically finite
and $f_1,\ldots,f_n$ define isolated base points 
in $\Proj(A)$. We then relate the definition ideal of $S$ to the
symmetric and the Rees algebras of 
the ideal $I=(f_1,\ldots,f_n) \subset A$, and prove some new acyclicity
criteria for the associated approximation complexes. Finally, we use
these results to obtain the implicit equation of $S$ in case $S$ is a
hypersurface, $\Proj(A)=\PP^{n-2}_k$ with $k$ a field, and base points
are either absent or local complete intersection isolated points.
\end{abstract}
\end{frontmatter}

\section{Introduction}
 Let $f_1,\ldots,f_{n}$ be non-zero homogeneous polynomials in
$k[X_1,\ldots,X_{n-1}]$ of the same degree $d$, where $k$ is a field. They
define a rational map
\begin{eqnarray*}
\PP^{n-2}_k & \xrightarrow{\lambda} & \PP^{n-1}_k \\
(X_1:\cdots:X_{n-1})& \mapsto & (f_1:\cdots:f_{n})(X_1:\cdots:X_{n-1}),
\end{eqnarray*}
whose closed image (the algebraic closure of the image) is a
hypersurface of $\PP^{n-1}_k$, say $H$, providing $\lambda$ is generically
finite. The computation of the equation of $H$, up to a non-zero
constant multiple,  is known as the
\emph{implicitization problem}, and is the main purpose of this
paper. We examine 
here this problem, under suitable assumptions, by studying  the
blow-up algebras associated to $\lambda$. More precisely, we use the
acyclicity of the approximation complexes introduced in \cite{SiVa81}
to develop a new algorithm giving the implicit equation $H$, in case
the ideal $(f_1,\ldots,f_{n})$ is zero-dimensional and locally a
complete intersection. This  
algorithm, based on determinants of complexes (also called MacRae's
invariants) 
computations, is the main contribution of this paper among other related results.

Implicitization in low dimensions, that is curve and surface
implicitizations, have significant applications to the field of computer aided
geometric design. For instance it helps drawing a curve or a surface nearby a
singularity, computing autointersection of offsets and drafts, or
 computing the intersection with parameterized curves or
 surfaces. Consequently 
the implicitization problem has been widely studied and is always an
active research area. It is basically an
elimination problem, and thus can be solved by
methods based on Gr\"obner basis computations, using for instance the
Buchberger's algorithm, and this  
in all generality. However, in practice it appears that such methods 
involve heavy computations and high complexity, and are hence rarely
used. For effective computations, algorithms from the old-fashioned
approach of resultants, going back to the end of the
nineteenth century, are preferred, even if
they apply only under particular conditions. Moreover they usually
give the implicit equation as a determinant of a matrix, which is a
very 
useful representation with well-known adapted
algorithms for numerical applications. We now survey these methods for
 curve and surface implicitizations in order to clarify what is 
known today, and also to shed light on the contribution of this paper.

\subsection{Overview on curve implicitization}\label{overcurve}
Curve implicitization is a quite easy problem which is well understood 
and solved. Let $\lambda$ be the rational map
\begin{eqnarray*}
\PP^1_k & \longrightarrow & \PP^2_k=\Proj(k[T_1,T_2,T_3])\\ 
(X_1:X_2)& \mapsto & (f_1:f_2:f_3)(X_1,X_2),
\end{eqnarray*}
where $f_1,f_2,f_3 \in k[X_1,X_2]$ are homogeneous polynomials of same degree $d\geq
1$ such that $\gcd(f_1,f_2,f_3)=1$, which implies that $\lambda$ is
regular (i.e.~there is no base points). Notice that this last hypothesis
is not restrictive because we can always obtain it by dividing each 
$f_i$, for $i=1,2,3$, by $\gcd(f_1,f_2,f_3)$. In fact, 
as we will see, the difficulty of the implicitization problem comes 
from the presence of base points, and the case of curve is
particularly simple since base points can be easily erased via a
single gcd
computation. There is roughly two types of methods to compute the
degree $d$ implicit equation of the closed image of
$\lambda$, which is of the form $C(T_1,T_2,T_3)^{\deg(\lambda)}$,
where $C$ is an irreducible homogeneous polynomial in
$k[T_1,T_2,T_3]$.

The first one is 
based on a resultant computation. We have the equality  
$$\Res(f_1-T_1f_3,f_2-T_2f_3)=C(T_1,T_2,1)^{\deg(\lambda)},$$
where $\Res$ denotes the classical resultant of two homogeneous
polynomials in $\PP^1_k$. 
We can thus obtain the implicit curve as the determinant of 
the well-known Sylvester's matrix of both polynomials $f_1-T_1f_3$ and 
$f_2-T_2f_3$ in variables $X_1$ and $X_2$. Another matricial
formulation is known to compute such a resultant, 
the Bezout's matrix which is smaller than the Sylvester's matrix. If
$P(X_1,X_2)$ and $Q(X_1,X_2)$ are 
two homogeneous polynomials in $k[X_1,X_2]$ of
the same degree $d$, the Bezout's matrix $\Bez(P,Q)$ is the matrix
$(c_{i,j})_{0\leq i,j \leq d-1}$ where the $c_{i,j}$'s are the coefficients of the
decomposition 
 $$\frac{P(S,1)Q(T,1)-P(T,1)Q(S,1)}{S-T}= \sum_{0\leq i,j \leq
   n-1}c_{i,j}S^iT^j.$$
Since $\det(\Bez(P,Q))=\Res(P,Q)$ (see e.g. \cite{Jou97}) we have 
$$\det(\Bez(f_1-T_1f_3,f_2-T_2f_3))=C(T_1,T_2,1)^{\deg(\lambda)}.$$
A variation of this 
formula, due to Kravitsky (see \cite{LKMV95}, theorem 9.1.11),  
is 
  $$\det(T_1\Bez(f_2,f_3)+T_2\Bez(f_3,f_1)+T_3\Bez(f_1,f_2))=C(T_1,T_2,T_3)^{\deg(\lambda)}$$
(see the appendix \ref{app} at the end of the paper for a proof in the general
context of  
anisotropic resultant).  
Notice that these resultant-based formulations have a big advantage, they give a
\emph{generic} implicitization 
formula, that is computations can be made with generic coefficients of 
polynomials $f_1,f_2,f_3$,  
and then all possible implicit equations are obtained as a
specialization of a generic formula.  
  
More recently a new method has been developed for curve
implicitization, called the ``moving lines'' method. It was introduced by
T. Sederberg and F. Chen in \cite{SeCh95}. This method is based on the
study of the first syzygies of the homogeneous ideal $I=(f_1,f_2,f_3)\subset k[X_1,X_2]$ 
in appropriate degrees (recall that we always suppose that
$\gcd(f_1,f_2,f_3)=1$). A moving line of degree $\nu$ is a homogeneous 
polynomial $$a_1(X_1,X_2)T_1+a_2(X_1,X_2)T_2+a_3(X_1,X_2)T_3,$$
where $a_1,a_2,a_3$ are homogeneous polynomials in $k[X_1,X_2]_\nu$,
which is said to follow the implicit  
curve if it vanishes when substituting each $T_i$ by $f_i$, 
 for $i=1,2,3$. In other words the triple $(a_1,a_2,a_3)$ is a syzygy
 of degree $\nu$ of the 
ideal $I$. Choosing $d$ moving lines $L_1,\ldots,L_d$ of degree $d-1$ following the
implicit curve,  we can construct a square matrix $M$ of size $d\times d$ corresponding to the
$k[T_1,T_2,T_3]$-module map
\begin{eqnarray*}
\oplus_{i=1}^{d}k[T_1,T_2,T_3] & \longrightarrow &
(k[X_1,X_2]_{d-1})\otimes_k k[T_1,T_2,T_3] \\
(p_1,p_2,\ldots,p_d) & \mapsto & \sum_{i=1}^d p_iL_i.
\end{eqnarray*}
It can be shown that one can always choose $L_1,\ldots,L_d$ such that the
determinant of $M$ is non-zero, and in this case it equals what one 
wants, that is $C(T_1,T_2,T_3)^{\deg(\lambda)}$ (see
\cite{Cox01}). This method hence  
yields an alternative to resultant computations, but however, it is
not a generic method.

For  
curve implicitization the two methods we have recalled   
work well and give the desired result, but, as we are now going to
see, it is no longer true for surface implicitization where the situation
becomes more intricate.

\subsection{Overview on surface implicitization}\label{oversurface}
Consider the rational map $\lambda$ 
\begin{eqnarray*}
\PP^2_k & \longrightarrow & \PP^3_k=\Proj(k[T_1,T_2,T_3,T_4])\\ 
(X_1:X_2:X_3)& \mapsto &
(f_1:f_2:f_3:f_4)(X_1,X_2,X_3), 
\end{eqnarray*}
where $f_1,f_2,f_3,f_4 \in k[X_1,X_2,X_3]$ are homogeneous polynomials
of same degree $d\geq 
1$, that we always suppose to be generically finite. Surface
implicitization problem  
consists in computing the implicit equation of the closed
image of $\lambda$, which is a hypersurface of $\PP^3_k$ of the form
$S(T_1,T_2,T_3,T_4)^{\deg(\lambda)}$, with $S$ irreducible. It is much 
more difficult than curve implicitization problem and, as far as we know, 
is unsolved in all generality (always without using Gr\"obner basis). 
Here again we suppose that 
$\gcd(f_1,f_2,f_3,f_4)=1$,  which is not restrictive as we have seen
in the preceding paragraph. However  
this hypothesis only erases base points 
in purely codimension 1 here, and hence it remains a finite set of
isolated base 
points, possibly empty, that we denote by $\BB$. We know then that the
equation of the closed image of $\lambda$ is a hypersurface of degree 
$(d^2-\sum_{p\in\BB}e_p)/\deg(\lambda)$, where $e_p$ denotes the
algebraic multiplicity of the point $p\in\BB$ (see theorem \ref{degree}, or
\cite{Cox01}, appendix). As for curve implicitization, two types of
methods have been developed to 
solve the surface implicitization problem: methods based on resultant
computations, and the method called ``moving surfaces''. We begin with
the first, which is also the oldest one.

If there is no base point (i.e.~$\BB=\emptyset$), one easily shows that 
$$\Res(f_1-T_1f_4,f_2-T_2f_4,f_3-T_3f_4)=S(T_1,T_2,T_3,1)^{\deg(\lambda)},$$
where $\Res$ denotes the classical resultant of three homogeneous
polynomials in $\PP^2$. This resultant can be computed with the well-known
Macaulay's matrices, but this involves gcd
computations since the determinant of each Macaulay's matrix gives
only a multiple of this resultant. In order to be efficient, direct
methods, avoiding these gcd computations, have been proposed. The general
case of an implicit surface without base points is solved in   
\cite{Jou96}, 5.3.17, using the anisotropic resultant; the implicit
hypersurface of $\lambda$ is then obtained as the determinant of  
a square matrix (notice that this
result was rediscovered in \cite{ArSe01}, in a more geometrical
setting, but only in some particular cases). As always with
resultant-based methods, the 
formula is \emph{generic}, i.e.~computations can be done with
indeterminate  
coefficients, giving the desired implicit equation as a specialization of a
generic equation. If base points exist the preceding resultant
vanishes identically, but we can use other constructions of resultants, taking
into account base points, as the reduced resultant (see \cite{Jou96},
5.3.29) or the residual resultant (see \cite{BusPhD}). In \cite{Bus01},
the residual resultant is used to solve the surface implicitization
problem if the 
ideal of base points is a local complete intersection satisfying 
some more technical conditions (see \cite{Bus01}, theorem 3.2; for 
instance, if the ideal $(f_1,f_2,f_3,f_4)$ is saturated the method
failed). This method has the disadvantage of requiring the knowledge of
the ideal of base points, but however it has the advantage of giving a
generic formula for all parameterizations admitting the same ideal of
base points.  

We now briefly overview the second announced method, the ``moving surfaces''
me\-thod which has also been introduced in \cite{SeCh95}. We first suppose that $\BB=\emptyset$. The two main         
ingredients of this method are 
moving planes and  moving quadrics. Like curves, a moving plane of
degree $\nu$ 
following the surface is a polynomial
in $k[X_1,X_2,X_3][T_1,T_2,T_3,T_4]$ of the form
$$a_1(X_1,X_2,X_3)T_1+a_2(X_1,X_2,X_3)T_2+a_3(X_1,X_2,X_3)T_3+a_4(X_1,X_2,X_3)T_4,$$
where $a_1,a_2,a_3$ and $a_4$ are homogeneous polynomials in
$k[X_1,X_2,X_3]_{\nu}$, such that it vanishes if we 
replace each $T_i$ by $f_i$. In a similar way, a moving quadric of degree $\nu$
following the surface is a polynomial
 of the form
$$a_{1,1}(X_1,X_2,X_3)T_1^2+a_{1,2}(X_1,X_2,X_3)T_1T_2+\ldots+a_{4,4}(X_1,X_2,X_3)T_4^2,$$
where $a_{i,j}$'s, with $1\leq i \leq j \leq 4$, are also homogeneous
polynomials in $k[X_1,X_2,X_3]_{\nu}$, such that it vanishes if we 
replace each $T_i$ by $f_i$. Choosing $d$ moving planes
$L_1,\ldots,L_n$ and $l=(d^2-d)/2$ moving
quadrics $Q_1,\ldots,Q_{l}$ of degree $d-1$ which follow the
surface, we can construct a square matrix, say $M$, corresponding to the
$k[\Tg]$-module map (where $k[\Tg]$ denotes $k[T_1,T_2,T_3,T_4]$)
\begin{eqnarray*}
\left(\oplus_{i=1}^{d} k[\Tg]\right) \bigoplus \left(\oplus_{j=1}^l k[\Tg]\right) & \rightarrow &
({k[X_1,X_2,X_3]}_{d-1})\otimes_k k[\Tg] \\ \nonumber
(p_1,\ldots,p_d,q_1,\ldots,q_l) & \rightarrow & \sum_{i=1}^d p_iL_i +
\sum_{j=1}^l q_iQ_i.\\  \nonumber
\end{eqnarray*} 
It can be shown that it is always possible to choose $L_1,\ldots,L_n$ and
$Q_1,\ldots,Q_{l}$ so that $\det(M)$ is non-zero and then equals 
$S^{\deg(\lambda)}$ (see \cite{Dan01}). In the recent
paper \cite{BCD02}, this method was improved: a similar matrix 
is constructed such that its 
determinant is the implicit surface, providing $\lambda$ is
birational and the base points ideal is locally a complete
intersection satisfying also some other technical conditions (for
instance the method rarely works if the ideal $(f_1,\ldots,f_4)$ is
saturated).   

\subsection{Contents of the paper} As we have already mentioned, this
article is centered on the implicitization problem. However we also
develop some other results which are closely related. After this quite 
long introduction, the paper is organized as follows.

Section 2 regroups some results on the geometry of the closed image of 
a rational map $\lambda$ from $\Proj(A)$, where $A$ is a $\ZZ$-graded
algebra, to $\Proj(k[T_1,\ldots,T_n])$, where $k$ is only supposed
to be a commutative ring, given by $n$ homogeneous elements
$f_1,\ldots,f_n$ in $A$ of the same degree $d\geq 1$. We first give
the ideal of definition of the scheme-theoretic image of $\lambda$,
that we will also called hereafter the closed image of $\lambda$, and
then prove that $\lambda$ can be extended to a projective morphism by
blowing-up. We end this section with a theorem on the degree of the
closed image of $\lambda$, giving a closed formula to compute it,
providing $k$ is a 
field, $A$ is a $\NN$-graded finitely generated $k$-algebra without
zero divisors, $\lambda$ is 
generically finite (onto its image) and its base points are isolated.

In section 3, we link the definition ideal of the closed image of
$\lambda$ to the symmetric and Rees algebras 
associated to the ideal $I=(f_1,\ldots,f_n)$ of the $\ZZ$-graded algebra 
$A$. We also investigate what happens when the ideal $I$ is of linear
type, that is the Rees algebra and the symmetric algebra of $I$ are
isomorphic. We roughly prove that it is possible to ``read'' the definition
ideal of the closed image of $\lambda$ from these blow-up algebras.  

The section 4 is devoted to the so-called approximation complexes
introduced by A. Simis and W.V. Vasconcelos in \cite{SiVa81}. In a
first part we introduce these complexes, recalling their definitions
and their basic, but very useful, properties. In a second part, having 
in mind to apply them to the implicitization problem, we
prove two new acyclicity criterions for the first approximation
complex associated to the ideal $I$, under suitable assumptions. The
last part of this section deals with a problem introduced by D. Cox and
solved by himself and H. Schenck  in
\cite{CoSc01}:  being given a homogeneous ideal
$I=(a_1,a_2,a_3)$ of $\PP^2$ of codimension two, the module of 
syzygies of $I$ vanishing at the scheme locus 
$V(I)$ is generated by the Koszul sygyzies if and only if $V(I)$ is a
local complete intersection in $\PP^2$. We show that this problem
appears naturally in 
the setting of the approximation complexes. This allows us to generalize
this result to higher dimensions, and observe that the ``natural''
condition on $I$ is not to be a local complete intersection, but to be
syzygetic. 

Finally, in section 5, we investigate the implicitization problem with
the tools we have developed in the previous sections. We suppose here
that $k$ is a field, $A$ is the polynomial ring $k[X_1,\ldots,X_{n-1}]$, and
$\lambda$ is generically finite (so that the closed image of $\lambda$ 
is a hypersurface of $\PP^{n-1}_k$). Using the
acyclicity criterions we proved in section 4, we show that an equation
of the closed image of $\lambda$ is obtained as the determinant of
each degree $\nu$, with $\nu\geq (n-2)(d-1)$, part of the
first approximation complex associated to $I$, under the hypothesis
that $I$ is either such that $V(I)=\emptyset$ or a local complete
intersection in $\Proj(A)$ of codimension $n-2$ (isolated base points). This new algorithm
 improves in particular the known techniques (not involving Gr\"obner basis) to
solve curve and surface implicitizations. For curve implicitization,
this algorithm is exactly the moving lines method in case
$\gcd(f_1,f_2,f_3)=1$, but its formalism shows clearly (in our point of
view) why $d-1$ is the 
good degree of moving lines to consider.  Moreover if we do not suppose $\gcd(f_1,f_2,f_3)=1$, then our
algorithm also applies. For surface implicitization, our algorithm
works if the ideal $I=(f_1,\ldots,f_4)$ defines local complete
intersection base points of codimension 2. Comparing with the method
of moving surfaces (and its improvement) and the resultant-based
methods, it is more general since both preceding methods sometimes
failed in this situation. We would like to mention that we have implemented 
our algorithm in the MAGMA language. This choice was guided by its 
needs: it just uses linear algebra, and more 
precisely linear systems solving routines.

\vspace{.2cm}

This work is based on a course given by Jean-Pierre Jouanolou at the
University Louis Pasteur of Strasbourg (France) during the academic
year 2000-2001.

\medskip

\noindent \textbf{Notation:} Let $R$ be a commutative ring. If $I$
is an ideal of $R$, we will denote $I^\thicksim$ its associated sheaf
of ideals on $\Spec(R)$. If the ring
$R$ is supposed to be $\ZZ$-graded, we will 
denote $\Proj(R)$ the quotient of $\Spec(R)\setminus V(\mm)$ by the 
multiplicative group $\GG_m$, where $\mm$ is the ideal of $R$
generated by all the homogeneous elements of non-zero degree. If $I$ 
is supposed to be an homogeneous ideal of $R$ we will denote $I^\#$
its associated sheaf of ideals on $\Proj(R)$. If in addition the ring $R$ is
supposed to be $\ZZ\times \ZZ$-graded, we will denote $\Biproj(R)$
the quotient of $\Spec(R)\setminus V(\tau)$ by the group
$\GG_m\times \GG_m$, where $\tau$ is the ideal of $R$ generated by all 
the products of two bihomogeneous elements of bidegree $(0,a)$ and $(b,0)$ 
such that $a\neq 0$ and $b\neq 0$. If $I$ is supposed to be a bihomogeneous
ideal of $R$, then we will denote $I^{\#\#}$ its associated sheaf of
ideals on $\Biproj(R)$.\\  
Finally, if $R$ is a commutative ring 
and if $J$ and $I$ are two ideals of $R$, we define the ideal of
inertia forms (or Tr\"agheitsformen) of $I$ with respect to $J$ to be  
$$\TF_J(I):=\cup_{n\in \NN}(I:_RJ^n)=\{a\in R \ \mathrm{such} \,
\mathrm{that} \ \exists 
n\in \NN  \  \forall \xi \in J^n \, : \,  a\xi \in I \}.$$

\section{The geometry of the closed image of a rational map}
Let $k$ be a commutative ring, $A$ a $\ZZ$-graded  
$k$-algebra, and denote by $\tau:k \rightarrow A_0$ the canonical
morphism of rings. Suppose given two integers $n\geq 1$ and $d\geq 1$,
and also an element $f_i \in A_d$ for each $i\in \{1,2,\ldots,n\}$.
The $k$-algebra morphism 
$$\begin{array}{cccc}
h : & k[T_1,\ldots,T_n] & \longrightarrow  & A \\
& T_i & \mapsto & f_i
\end{array}$$
yields a $k$-scheme morphism 
$$\begin{array}{cccc}
\mu : & \cup_{i=1}^n D(f_i) & \longrightarrow  & \cup_{i=1}^n
D(T_i)=(\AA_k^n)^*, \end{array}$$
where $\AA_k^n$ denotes the affine space of dimension $n$ over $k$,
and $(\AA_k^n)^*$ this affine space without the origin.
If we quotient the morphism $\mu$ by the action of the multiplicative group
$\GG_m$, we obtain another $k$-scheme morphism, 
\begin{equation}\label{deflambda}
\lambda \ : \  \cup_{i=1}^n D_+(f_i)  \longrightarrow
\PP_k^{n-1}=\Proj(k[T_1,\ldots,T_n]), \end{equation}
and the following commutative diagram:
$$\begin{CD}
D(\fg):=\cup_{i=1}^n D(f_i) @>\mu>> (\AA_k^n)^* \\
@V{\Proj}VV  @V{\Proj}VV \\
D_+(\fg):=\cup_{i=1}^n D_+(f_i) @>\lambda>> \PP_k^{n-1}. 
\end{CD}$$

The main aim of this paper is to provide a ``formula'' for the closure
of the set-theoretic images of the maps $\lambda$ and $\mu$ in case
they are hypersurfaces. It turns
out that, under suitable assumption, the closure of a set-theoretic
image of a scheme morphism have a natural
scheme structure called the  
\emph{scheme-theoretic image}. More precisely, if $X\xrightarrow{\rho} Y$ is a
morphism of schemes, then its scheme-theoretic image exists if, for
instance, $\rho_*(\OO_X)$ is quasi-coherent. It is the smallest closed 
subscheme of $Y$ containing the set-theoretic image of $\rho$. Its
associated ideal sheaf is given by 
$\ker(\rho^\#:\OO_Y\rightarrow \rho_*\OO_X)$ and it is supported on the
 closure of its set-theoretic image (see \cite{GrDi71}
 I.6.10.5 for more details). The first part of this section is    
devoted to the description of the scheme-theoretic image of the map 
$\lambda$ (and $\mu$): we give its sheaf of ideals on $\PP^{n-1}_k$, 
and a projective morphism which extends it. In a second part we
prove a formula to compute the degree of the    
scheme-theoretic image of $\lambda$ in case elements 
$f_1,\ldots,f_n$ define a finite number of points in $A$ (possibly
zero). This formula generalizes the 
one commonly uses when dealing with the implicitization problem to
obtain the degree of a parameterized surface of 
$\PP^3$ (see e.g. \cite{Cox01}, appendix). 

\subsection{The ideal of the closed image of $\lambda$} We begin here by the
description of the ideal sheaf of the scheme-theoretic image, that we
will also often called the closed image, of $\lambda$ (and $\mu$). 

\begin{theorem} The scheme-theoretic images of the respective morphisms $\mu$ and  $\lambda$ are given by
   $$V(\ker(h)^\thicksim \,
  )_{|(\AA_k^n)^*}=V\left(\TF_{(T_1,\ldots,T_n)}(\ker (h))^{\thicksim}
    \,\right)_{|(\AA_k^n)^*},$$
and
$$V(\ker(h)^\# \,)=V\left(\TF_{(T_1,\ldots,T_n)}(\ker (h))^\#\,\right).$$
\end{theorem}  
\begin{proof} First notice that for all $i\in \{1,2,\ldots,n\}$, the
  restriction 
$\mu_{|D(f_i)}$ is the spectrum of the morphism of rings
$ k[T_1,\ldots,T_n]_{T_i} \rightarrow A_{f_i}$,
obtained from $h$ by localization. Since 
$\mu^{-1}(D(T_i))=D(f_i)$ for all $i\in\{1,2,\ldots,n\}$, it follows that
$\mu$ is an affine morphism.
Now denote by $j$ the canonical inclusion $(\AA_k^n)^* \hookrightarrow
\AA_k^n$. The sheaf $\OO_{D(\fg)}$ is quasi-coherent and the morphism
$j\circ \mu$ is quasi-compact and separated (since it is affine), so we
 deduce that $(j\circ
\mu)_*(\OO_{D(\fg)})$ is quasi-coherent (see \cite{Har77},
prop. II.5.8), and we obtain
$$(j\circ\mu)_*(\OO_{D(\fg)})=\ker( \prod_{i=1}^n A_{f_i}
  \rightarrow \prod_{1\leq i,j \leq n, \, i\neq j} A_{f_if_j})^\thicksim,$$
and also
$$\mu_*(\OO_{D(\fg)})=\ker( \prod_{i=1}^n A_{f_i}
  \rightarrow \prod_{1\leq i,j \leq n, \, i\neq j}
  A_{f_if_j}{)^\thicksim}_{|(\AA_k^n)^*}.$$
  Notice that $\ker( \prod A_{f_i}
  \rightarrow \prod  A_{f_if_j})$ is a $k[T_1,\ldots,T_n]$-module via the map $h$. 
From here the kernel of the canonical morphism $\OO_{(\AA_k^n)^*}
\rightarrow \mu_*\OO_{D(\fg)}$ is exactly
\begin{eqnarray*}
&  & \ker ( \ k[T_1,\ldots,T_n]\xrightarrow{h} \ker( \prod_{i=1}^n A_{f_i}
  \rightarrow \prod_{1\leq i,j \leq n}
  A_{f_if_j}) \ {)^\thicksim}_{|(\AA_k^n)^*} \\
& = & \ker(k[T_1,\ldots,T_n]\xrightarrow{h}
\frac{A}{H^0_{(f_1,\ldots,f_n)}(A)}{)^\thicksim}_{|(\AA_k^n)^*}\\
& = & (\TF_{(T_1,\ldots,T_n)}(\ker (h)){)^\thicksim}_{|(\AA_k^n)^*}\\
& = & {\ker(h)^\thicksim}_{|(\AA_k^n)^*},
\end{eqnarray*}
and this proves the first statement of the proposition. If we quotient
by the action of the multiplicative group $\GG_m$ we obtain
$$\ker(\OO_{\PP_k^{n-1}} \rightarrow
\lambda_*(\OO_{D_+(\fg)}))=(\TF_{(T_1,\ldots,T_n)}(\ker
(h)))^\#=\ker(h)^\#$$
which is the second statement.
\end{proof}

\begin{remark}
The ideal $\TF_{(T_1,\ldots,T_n)}(\ker(h))$ is often called the
saturated ideal of $\ker(h)$ in $k[T_1,\ldots,T_n]$. Also, in
the case where $H^0_{(f_1,\ldots,f_n)}(A)=0$ we have 
$\ker(h)=\TF_{(T_1,\ldots,T_n)}(\ker(h))$, i.e.~$\ker (h)$ is
saturated, and hence, for all $\nu \in \ZZ$, 
$\ker(h)_\nu \simeq H^0(\PP_k^{n-1},\ker(h)^\#(\nu))$.
\end{remark}

Now we show that the map $\lambda$ can be extended
to a projective morphism whose image is the closed image of $\lambda$. This
result is quite classical (see for instance \cite{Har77}, II.7.17.3),
but we prefer to detail it, giving a particular attention to
graduations.  

First of all we set the indeterminates $T_1,\ldots,T_n$ to degree 1,
the usual choice to construct
$\PP^{n-1}_k=\Proj(k[T_1,\ldots,T_n])$. It follows that the ring
$$A[T_1,\ldots,T_n]=A\otimes_k k[T_1,\ldots,T_n]$$ is naturally 
bi-graded by the tensor product graduation associated to the
graduations of $A$ and $k[T_1,\ldots,T_n]$. Let $Z$ be a new
indeterminate and set its degree to $1-d$. In this way the
$\ZZ\times \ZZ$-graded $A$-algebra morphism
\begin{eqnarray}\label{Rees}
A[T_1,\ldots,T_n] & \longrightarrow & A[Z,Z^{-1}] \\
T_i & \mapsto & f_iZ \nonumber
\end{eqnarray}
is homogeneous, and hence gives a natural
$\ZZ\times\ZZ$-graduation to its image: the Rees algebra
$\Rees_A(I)$, where $I=(f_1,\ldots,f_n)\subset A$. For this graduation, we
denote by $Bl_I$ the blow-up
$$Bl_I=\mathrm{Biproj}(\Rees_A(I)) \subset
\mathrm{Biproj}(A[T_1,\ldots,T_n])=\Proj(A)\times_k\PP_k^{n-1},$$ and
respectively by $u$ and $v$ the two canonical projections
$$Bl_I\xrightarrow{u} \Proj(A) \ \ \mathrm{and} \ \
Bl_I\xrightarrow{v} \PP_k^{n-1}.$$
Since $I_p=A_p$ for all $p \in D_+(\fg)$, the
restriction of $u$ to $\Omega=u^{-1}(D_+(\fg))$ is an isomorphism from
$\Omega$ to $D_+(\fg)$. The following proposition shows that the
morphism of $k$-schemes $v$ extends the morphism $\lambda$ by
identifying $\Omega$ with $D_+(\fg)$ via $u_{|\Omega}$.

\begin{proposition}\label{commute}
With the above notations, the following diagram is commutative:
\begin{diagram}
\Omega & & \\
\dTo^{u_{|\Omega}}_{\wr} & \rdTo^{v_{|\Omega}} & \\
D_+(\fg)& \rTo{\lambda} & \PP_k^{n-1}
\end{diagram}
\end{proposition}
\begin{proof}
Setting the weights of the indeterminates $T_1,\ldots,T_n$ to $d$ and
$Z$ to 0, we obtain two graded $k$-algebra morphisms
\begin{eqnarray*}
incl \circ h : k[T_1,\ldots,T_n] & \xrightarrow{h} A \hookrightarrow &
\Rees_A(I) \subset A[Z,Z^{-1}]\\
T_i & \mapsto & f_i
\end{eqnarray*}
and
\begin{eqnarray*}
\rho : k[T_1,\ldots,T_n] & \longrightarrow & \Rees_A(I) \subset A[Z,Z^{-1}]\\
T_i & \mapsto & f_iZ.
\end{eqnarray*}
Let $i\in \{1,2,\ldots,n\}$ be fixed. Over $D_+(T_i)\subset
\PP^{n-1}_k$, the map $\lambda\circ u_{|\Omega}$ is affine, associated
to the $k$-algebra morphism
\begin{eqnarray*}
(incl\circ h)_i : k[T_1,\ldots,T_n]_{(T_i)} & \longrightarrow &
(\Rees_A(I))_{(f_i)} \hookrightarrow A_{(f_i)}[Z,Z^{-1}]
\\
T_j/T_i & \mapsto & f_j/f_i,
\end{eqnarray*}
and the map $v_{|\Omega}$ is also affine, associated to the
$k$-algebra morphism
\begin{eqnarray*}
\rho_i : k[T_1,\ldots,T_n]_{(T_i)} & \longrightarrow & \Rees_A(I)_{(f_iZ)} \subset A_{(f_i)}[Z,Z^{-1}]\\
T_j/T_i & \mapsto & f_jZ/f_iZ.
\end{eqnarray*}
This clearly shows that $\rho_i(T_j/T_i)=f_j/f_i=(incl\circ h)_i$ in
$A_{(f_i)}[Z,Z^{-1}]$, which implies the proposition.
\end{proof}


\subsection{The degree of the closed image of $\lambda$}
In what follows we give an explicit formula to compute the
degree of the closed image of $\lambda$, that is the image of $v$,
providing that $\Proj(A/I)$ is a zero-dimensional scheme (possibly
empty), where $I$ 
denotes the ideal $(f_1,\ldots,f_n)$. For this we
first fix some 
notation about multiplicities (we refer to \cite{BrHe93} for
a complete treatment on the subject), and recall some classical
results. 

From now we suppose that $k$ is a field and we denote by $C$ the
$\NN$-graded polynomial ring $k[X_1,\ldots,X_r]$ with $\deg(X_i)=1$
for all $i=1,\ldots,r$. For all
$\ZZ$-graded finite $C$-module $M$ one defines the \emph{Hilbert series} of $M$: 
$$\H_M(T)=\sum_{\nu\in\ZZ} \dim_k(M_\nu)T^\nu.$$
If $\delta$ denotes the Krull dimension of $M$, there exists a unique
polynomial $L_M(T)$ such that $L_M(1)\neq 0$ and
$$\H_M(T)=\sum_{\nu\in\ZZ}
\dim_k(M_\nu)T^\nu=\frac{L_M(T)}{(1-T)^\delta}.$$
The number $L_M(1)$ is an invariant of the module $M$ called the
\emph{multiplicity} of $M$; we will denote it by $\mult_k(M):=L_M(1)$. Another
way to obtain this invariant is the $\emph{Hilbert polynomial}$ of
$M$, denoted $\P_M(X)$. It is a polynomial of degree $\delta-1$ 
such that $\P_M(\nu)=\dim_k(M_\nu)$ for all sufficiently large
$\nu\in\NN$. The Hilbert polynomial is of the form
$$\P_M(X)=\frac{a_{\delta-1}}{(\delta-1)!}X^{\delta-1}+\ldots+a_0,$$
and we have the equality $\mult_k(M)=L_M(1)=a_{\delta-1}$. One 
can also define it as the Euler characteristic of $M^\sharp$ on 
$\PP_k^{r-1}$, that is for all $\nu\in \ZZ$ we have:
\begin{equation}\label{ki}
\P_M(\nu)=\chi(\PP^{r-1}_k,M^\sharp(\nu))=\sum_{i\geq 0}(-1)^i
\dim_k H^i(\PP^{r-1}_k,M^\sharp(\nu)).
\end{equation}
A well-known formula relates the Hilbert
series and the Hilbert polynomial: for all $\nu\in \ZZ$ 
\begin{equation}\label{hilb}
 H_M(T)_{|T^\nu}-P_M(\nu)=\sum_{i\geq 0} (-1)^i \dim_k
H^i_\mm(M)_\nu,
\end{equation}
where $H_M(T)_{|T^\nu}$ is the coefficient of $T^\nu$ in the Hilbert
series $H_M(T)$, that is $\dim_k(M_\nu)$, and $\mm$ is the
irrelevant ideal of $k[X_1,\ldots,X_r]$, that is
$\mm=(X_1,\ldots,X_r)$. Recall also that 
for all $i>0$ 
$$H^i(\PP^{n-1}_k,M^\sharp(\nu))\simeq H^{i+1}_\mm(M)_\nu,$$
and that, for all $\nu\in \ZZ$, we have an exact sequence:
\begin{equation}\label{seHm}
0\rightarrow H^0_\mm(M)_\nu\rightarrow M_\nu \rightarrow
H^0(\PP^{n-1}_k,M^\sharp(\nu)) \rightarrow H^1_\mm(M)_\nu \rightarrow
0.
\end{equation}

Such a definition of multiplicity for $M$ is called a \emph{geometric}
multiplicity and is also often called the
degree of $M$ because of its geometric meaning. Indeed let $J$ be a
graded ideal of $C$ and consider the quotient ring $R=C/J$. If $\delta$
denotes the dimension of $R$ then the subscheme $\Proj(R)$ of $\PP^{n-1}_k$
is of dimension $\delta-1$. The degree of $\Proj(R)$ over $\PP^{n-1}_k$ is
defined to be the number of points obtained by cutting $\Proj(R)$
by $\delta-1$ generic linear forms. To be more precise, if
$l_1,\ldots,l_{\delta-1}$ are generic linear forms of $\PP^{n-1}_k$, then
the scheme $S=\Proj(R/(l_1,\ldots,l_{\delta-1}))$ is finite and we set
$$\deg_{\PP^{r-1}_k}(\Proj(R)):=\dim_k
\Gamma(S,\OO_S)=\dim_k\left(\frac{R_\nu}{{(l_1,\ldots,l_{\delta-1})}_\nu}\right),$$
for all sufficiently large $\nu$.
This geometric degree is in fact exactly what
we called the multiplicity of $R$, i.e.~we have
$$\deg_{\PP^{r-1}_k}(\Proj(R))=\mult_k(R)=L_M(1)=a_{d-1}.$$ 
This equality follows immediately from the exact sequence
$$ 0 \rightarrow R(-1) \xrightarrow{\times l_1} R \rightarrow R/(l_1)
\rightarrow 0$$
which gives $\H_{R/(l_1)}(T)=(1-T)\H_R(T)=L_M(T)/(1-T)^{\delta-1}$, and an
easy recursion.

At this point we can compute the multiplicity of a
Veronese module which will be useful later.

\begin{lemma}\label{verodeg}
  Let $M$ be a $\ZZ$-graded finite $C$-module. 
We denote by $M^{(d)}$, for all integers $d$, the Veronese 
module $\bigoplus_{\nu \in \ZZ} M_{d\nu}$ which is a
$k[X_1^d,\ldots,X_n^d]$-module. Setting $\delta=\dim_k(M)$ we have:
$$\mult_k(M^{(d)})=d^{\delta-1}\mult_k(M).$$
\end{lemma}
\begin{proof}
  Let $d$ be a fixed integer. It is known that for all integers $\nu$
  the sum  $\sum_{\xi^d=1} \xi^\nu$, where the sum is taken over all
  the $d$-unit roots, equals $d$ if $\nu$ divides $d$
  and zero otherwise. Using this result one deduces
  $$\H_{M^{(d)}}(T^d)=\frac{L_{M^{(d)}}(T^d)}{(1-T^d)^\delta}=\frac{1}{d}\sum_{\xi^d=1} \H_M(\xi T)=\frac{1}{d}\sum_{\xi^d=1}\frac{L_M(\xi T)}{(1-\xi T)^\delta}.$$ 
We obtain that $L_{M^{(d)}}(T^d)$ equals 
$$\frac{1}{d}\sum_{\xi^d=1} \frac{(1-T^d)^\delta}{(1-\xi
  T)^\delta}L_M(\xi T)=\frac{1}{d}\sum_{\xi^d=1}
\frac{(1-T)^\delta(1+T+\ldots+T^{d-1})^\delta}{(1-\xi T)^\delta}L_M(\xi T).$$
As $\mult_k(M^{(d)})=L_{M^{(d)}}(1)$ it remains only to compare 
 the preceding equality when $T\mapsto 1$, and one obtains
$$\mult_k(M^{(d)})=\frac{d^\delta}{d}L_M(1)=d^{\delta-1}\mult_k(M),$$
the desired result.
\end{proof}

In order to state the main theorem of this section we need to
introduce another notion of multiplicity which is called the
\emph{algebraic} multiplicity. Let $(R,\mm)$ be a local noetherian ring and 
$M\neq 0$ a finite $R$-module. Let $I\subset
\mm$ be an ideal of $R$ such that there exists an integer $t$ satisfying
$\mm^tM\subset IM$ (any such ideal is called a definition ideal of
$M$), the numerical function $\length(M/I^{\nu}M)$ is a polynomial
function for sufficiently large values of $\nu \in \NN$. This polynomial,
denoted $\mathrm{S}_M^I(X)$, is
called the \emph{Hilbert-Samuel} polynomial of $M$ with respect to
$I$. It is of degree $\delta=\dim_k(M)$ and of the form:
$$\mathrm{S}_M^I(X)=\frac{e(I,M)}{\delta!}X^\delta+\mathrm{ \it terms \ of \ lower \
  powers \ in \ } X.$$
The algebraic multiplicity of $I$ in $M$ is the number
$e(I,M)$ appearing in this polynomial. With such a definition of
algebraic multiplicity one can define the algebraic multiplicity of a
zero-dimensional scheme 
as follows: let $J$ be a graded ideal of a $\NN$-graded ring $R$, then if
$T=\Proj(R/J)$ is a finite subscheme of $\Proj(R)$, its algebraic
multiplicity is 
$$e(T,\Proj(R))=e(J^\sharp,R^\sharp)=\sum_{t\in
  T}e(J^\sharp_t,\OO_{\Proj(R),t})=\sum_{t\in T}e(J_t,R_t).$$
We are now ready to state the main result of this section.

\begin{theorem}\label{degree} Suppose that $k$ is a field and $A$ a
  $\NN$-graded $k$-algebra of the form $k[X_1,\ldots,X_r]/J$, where
  $J$ is a prime 
  homogeneous ideal and each $X_i$ is of degree one. Denote by
  $\delta$ the dimension of $A$ and let $I=(f_1,\ldots,f_n)$ be an
  ideal of $A$ such that each $f_i$ is of degree $d\geq 1$. Then,
  if $T=\Proj(A/I)$ is finite over $k$, the number 
  $d^{\delta -1}\deg_{\PP^{r-1}_k}(\Proj(A))-e(T,\Proj(A))$ equals
  $$\left\{\begin{array}{l} 
    \deg(\lambda).\deg_{\PP^{n-1}_k}(S) \ {\rm if} \ \lambda \ \mbox{{\rm is
        generically finite}}\\ 
    0 \ \mathrm{if} \ \lambda \ {\rm is \ not \ generically \ finite},
      \end{array} \right.$$
where $S$ denotes the closed image of $\lambda$ (i.e.~the image of $v$).
\end{theorem}
\begin{proof}
  We denote by $C$ the polynomial ring $k[X_1,\ldots,X_r]$ and
  consider it as a $\NN$-graded 
  ring by setting $\deg(X_i)=1$ for all $i=1,\ldots,r$. We denote also 
  by $\mm$ its irrelevant ideal $\mm=(X_1,\ldots,X_r)$, and for
  simplicity we set $X=\Proj(A)$. The
  $k$-algebra $A$ is a $\NN$-graded $C$-module and combining
  \eqref{ki} and \eqref{hilb} we obtain for all $\nu \in \NN$
  \begin{equation}\label{F-A}
    \dim_k(A_{d\nu})=\chi(X,\OO_X(d\nu))+\sum_{i\geq 0} (-1)^i
\dim_k H^i_\mm(A)_{d\nu}.
  \end{equation}
Similarly, for all $\nu \in \NN$ the ideal $I^\nu$ of $A$ is also a
$\NN$-graded $C$-module. We hence deduce that for all $\nu \in \NN$ 
  \begin{equation}\label{F-I}
    \dim_k((I^\nu)_{d\nu})=\chi(X,(I^\nu)^\sharp(d\nu))+\sum_{i\geq 0}
    (-1)^i  \dim_k H^i_\mm(I^\nu)_{d\nu}.
  \end{equation}
Now the exact sequence of sheaves on $X$
$$0\rightarrow (I^\nu)^\#(d\nu) \rightarrow \OO_X(d\nu)
\rightarrow \frac{\OO_X}{(I^\nu)^\#}(d\nu) \rightarrow 0,$$
shows that, always for all $\nu \in \NN$, 
$$\chi(X,\OO_X(d\nu)) -
\chi(X,(I^\nu)^\sharp(d\nu)) = 
\chi(X,\frac{\OO_X}{(I^\nu)^\#}(d\nu)).$$
Since $T=\Proj(A/I)$ is supposed to be finite, for all
 integers $\nu$ we have  
$$\chi(X,\frac{\OO_X}{(I^\nu)^\sharp}(d\nu))= \chi(X, \OO_X/(I^{\nu})^\sharp) =
\dim_k H^0(X,\OO_X/(I^{\nu})^\sharp),$$
where the last equality comes from \eqref{ki}.
Subtracting \eqref{F-I} to 
\eqref{F-A}, it follows that for all integers $\nu$
{\small \begin{eqnarray*}
\dim_k(A^{(d)}_\nu)-\dim_k((I^\nu)_{d\nu}) & = & \dim_k
H^0(X,\OO_X/(I^{\nu})^\sharp) + \sum_{i\geq 0} (-1)^i
\dim_k H^i_\mm(A)_{d\nu}\\
 & & - \sum_{i\geq 0}
    (-1)^i  \dim_k H^i_\mm(I^\nu)_{d\nu},
\end{eqnarray*}}

\noindent and for all $\nu \in \NN$ sufficiently large (recall that
$H^i_\mm(A)_\nu=0$ for $i\geq 0$ and $\nu \gg 0$)
{\small \begin{equation}\label{bigF}
\dim_k(A^{(d)}_\nu)-\dim_k((I^\nu)_{d\nu})  =  \dim_k
H^0(X,\OO_X/(I^{\nu})^\sharp)  - \sum_{i\geq 0}
    (-1)^i  \dim_k H^i_\mm(I^\nu)_{d\nu}.
\end{equation} }
At this point, notice that the algebraic multiplicity $e(T,X)$ appears
naturally in this formula. Indeed, for all $\nu \gg 0$,
$$\dim_k H^0(X,\OO_X/(I^{\nu})^\sharp)=\sum_{t\in T}\length(A_t /
I^\nu_t)=\frac{e(T,X)}{\delta!}\nu^\delta + R(\nu),$$
where $R$ is a polynomial in $k[\nu]$ of degree stricly less than
$\delta$.  

Now we relate \eqref{bigF} with the Rees algebra of $I$. 
Recall that we have a map $h:k[T_1,\ldots,T_n]\rightarrow A$ which
sends each $T_i$ to the polynomial $f_i\in A_d$. The polynomial ring
$k[T_1,\ldots,T_n]$ is naturally $\NN$-graded by setting $\deg(T_i)=1$
for $i=1,\ldots,n$. We have the following exact sequence of
$\NN$-graded $k[T_1,\ldots,T_n]$-modules
$$0\rightarrow \ker(h) \rightarrow k[T_1,\ldots,T_n] \xrightarrow{h}
\Im(h) \rightarrow 0.$$
From the definition of the map $h$ it comes $\Im(h)_\nu=I^\nu\cap
A_{d\nu}=(I^\nu)_{d\nu}$, and hence 
$\dim_k((I^\nu)_{d\nu})=\dim_k(\Im(h)_\nu)$, for all $\nu \in \NN$.
The Rees algebra of $I$, that we denote hereafter by $B$, is the image of
the following $A$-algebra 
morphism
\begin{eqnarray*}
  A[T_1,\ldots,T_n] & \xrightarrow{\beta} & B:=\Rees_I(A) \subset A[Z]\\
  T_i & \mapsto & f_iZ.
\end{eqnarray*}
As $A$ is $\NN$-graded, $B$ is naturally $\ZZ\times \ZZ$-graded if
$\beta$ is a homogeneous map, that is if $\deg(T_i)=d+\deg(Z)$, for all
$i=1,\ldots,n$. We choose $\deg(T_i)=0$ for $i=1,\ldots,n$, and hence
$\deg(Z)=-d$ (notice that $B$ has its own grading as a
$A$-algebra, so this choice do not interfere with the grading of
$C$). In this way $B_{p,q}=I^q_{p+dq}$ for all $(p,q)\in \ZZ
\times \ZZ$, and hence $B_{p,\bullet}=\bigoplus_{q\in \ZZ}
I^q(dq)_p$. It follows that $B=\bigoplus_{q\geq 0} I^q(dq)$
as a $\NN$-graded $A$-module. Hereafter we will denote by $B_p$ the
graded part of $B$ as a $A$-module, that is $B_p:=B_{p,\bullet}$. It
is then easy to check that $B_0=\bigoplus_{q\geq 0} (I^q)\cap
A_{dq}=\Im(h)$ and 
$\bigoplus_{q\geq 0} H^i_\mm(I^q)_{dq}=H^i_\mm(B)_0$, both as
$\NN$-graded $k[T_1,\ldots,T_n]$-modules. Finally \eqref{bigF} can be
rewritten as 
{\small \begin{equation}\label{NbigF}
\dim_k(A^{(d)}_\nu)-\dim_k H^0(X,\OO_X/(I^{\nu})^\sharp)=\dim_k(B_{0,\nu})-
\sum_{i\geq 0}(-1)^i\dim_k(H^i_\mm(B)_{0,\nu}),
\end{equation}}
always for all integers $\nu$ sufficiently large.

We are now ready to prove the second assertion of this
theorem. Denoting $\pp=(\ker(h))\subset k[T_1,\ldots,T_n]$, if we
suppose that 
the map $\lambda$ is not generically finite, then
$\dim(k[T_1,\ldots,T_n]/\pp)<\delta=\dim(A)$. Since $B_0$ and
$H^i_\mm(B)_0$, with $i\geq 0$, are $k[T_1,\ldots,T_n]/\pp$-modules,
we deduce that $B_{0}$ and $H^i_\mm(B)_{0}$, for all $i\geq 0$, have
Hilbert polynomials of 
degree strictly less than $\dim(A)$. But the left side of equality
\eqref{NbigF} is a polynomial of degree $\delta$ with leading
coefficient $$\frac{d^{\delta-1}\mult_k(A)-e(T,\Proj(A))}{\delta!},$$
and thus we obtain the second statement of this theorem (notice that
we did not use here the hypothesis $J$ is prime).

\medskip

We now concentrate on the first point of this theorem. We denote $\Proj_{{\rm
    A}}(B)$ the $\Proj$ of $B$ as a $\NN$-graded $A$-module, so that
\eqref{seHm} gives (in a more general setting, see \cite{Eis94} theorem A4.1), for all $\nu\in\NN$, 
the exact sequence of $\NN$-graded
$k[T_1,\ldots,T_n]$-modules:  
\begin{equation}\label{seB}
  0\rightarrow H^0_\mm(B)_{0} \rightarrow B_{0} \rightarrow
\Gamma(\Proj_{{\rm A}}(B),\OO_{\Proj_{{\rm A}}(B)}) \rightarrow H^1_\mm(B)_{0}
\rightarrow 0.
\end{equation}

Now remark that $\lambda$ is generically finite if and only if the 
morphism $\overline{v}:\Proj_A(B)\rightarrow
\Spec(k[T_1,\ldots,T_n]/\pp)$, associated to the canonical morphism
$k[T_1,\ldots,T_n]/\pp\rightarrow B$ sending each $T_i$ to $f_iZ$, is
also generically finite. Thus there 
exists an element $L \in k[T_1,\ldots,T_n]/\pp$ such that $\overline{v}_L$
is finite. We obtain that $H^i_\mm(B_L)=0$ for all $i\geq 2$, since
$\Proj_A(B)_{|D(L)}$ is affine and hence
$H^i(\Proj_A(B)_{|D(L)},{\OO_{\Proj_A(B)}}_{|D(L)})=0$ for all $i\geq
1$ (see \cite{Har77}, III.3), and it follows that, for all $i\geq 2$,
$H^i_\mm(B)_{0}$ has Hilbert polynomial of
degree less than or equal to $\delta-1$. This and \eqref{seB} show that
\eqref{NbigF} 
reduces to the equality        
 \begin{equation}\label{NbigFbis}
\dim_k(A^{(d)}_\nu)-\dim_k H^0(X,\OO_X/(I^{\nu})^\sharp)=\dim_k \Gamma(\Proj_{{\rm A}}(B),\OO_{\Proj_{{\rm A}}(B)})_\nu,
\end{equation}
which is true for all $\nu \gg 0$.

Consider now the generic point $s$ of the image $S=V(\pp^\#)\subset
\PP^{n-1}_k$ of $v$, and denote by $Y$ the fibred product 
$$Y:=Bl_{\PP^{n-1}_k}\times_{\PP^{n-1}_k}
\Spec(\OO_{\PP^{n-1}_k,s})=Bl_I\times_S \Spec(\OO_{S,s}).$$
We claim that, if the map $\lambda$ is generically finite, there exists an isomorphism of $\left(\frac{k[T_1,\ldots,T_n]}{\pp}\right)_\pp$-modules: 
\begin{equation}\label{iso}
\Gamma(\Proj_A(B),\OO_{\Proj_A(B)})_\pp\simeq
\left(\frac{k[T_1,\ldots,T_n]}{\pp}\right)_\pp\otimes_{\OO_{S,s}}
\Gamma(Y,\OO_{Y}).
\end{equation}
Suppose for the moment that this claim is true. From proposition 
\ref{commute} $\lambda$ is generically finite if and only if
$v$ is. Supposing now that $J$ is prime, this is equivalent to say
that the scheme $Y$ is the spectrum of a field 
which is a finite extension of the field $\OO_{S,s}$; the degree of
$\lambda$ is then the degree of the extension $\Gamma(Y,\OO_Y)$ of
$\OO_{S,s}$. By \eqref{iso} we deduce that there exists a morphism of
$\NN$-graded $k[T_1,\ldots,T_n]/\pp$-modules of finite type
$$\gamma : \frac{k[T_1,\ldots,T_n]}{\pp}\otimes_{\OO_{S,s}}\Gamma(Y,\OO_Y)
\rightarrow \Gamma(\Proj_A(B),\OO_{\Proj_A(B)}),$$
which becomes an isomorphism by localisation at $\pp$. The kernel $E$ and
the cokernel $F$ of $\gamma$ are hence graded
$k[T_1,\ldots,T_n]/\pp$-modules of finite type annihilated by an
homogeneous element of non-zero degree. It follows that $\dim(E)$ and
$\dim(F)$ are strictly lower than $\dim(A)$, which equals $\dim(S)$ by
hypothesis. From here we may deduce that the Hilbert polynomial of the
$k[T_1,\ldots,T_n]$-module $\Gamma(\Proj_A(B),\OO_{\Proj_A(B)})$ is of
degree $\delta$ with $\deg(\lambda)\deg(S)/\delta!$ as leading 
coefficient, since the source
of $\gamma$ consists in $\deg(\lambda)$ copies of
$k[T_1,\ldots,T_n]/\pp$ which has Hilbert polynomial of degree
$\delta$ with leading coefficient $\deg(S)/\delta!$. From \eqref{NbigFbis} 
the first statement of the theorem follows. To complete the proof it
thus remains to prove the claim \eqref{iso}. Hereafter, for
convenience,  we will write, as usually,   
$k[\Tg]$ for the polynomial ring $k[T_1,\ldots,T_n]$.

\medskip

We denote by $\hat{v} : \Proj_A(B) \rightarrow \AA^n_k=\Spec(k[\Tg])$
the canonical projection, whose closed image is $V(\pp^\thicksim)$,
which is $\{1\}\times \GG_m$-equivariant. Recall that, by definition,
$$Bl_I:=\Biproj(B)=(\Proj_A(B)_{|(\AA^n_k)^*})/\{1\}\times \GG_m,$$
and that the morphism $v:Bl_I\rightarrow \PP^{n-1}_k$ is obtained from
the morphism
$\hat{v}_{|(\AA^n_k)^*}:\Proj_A(B)_{|(\AA^n_k)^*}\rightarrow
(\AA^n_k)^*$ by passing to the quotient. In fact, as the variables
$T_i$ are all of degree 1, we have the following commutative diagram
\begin{diagram}
\Proj_A(B)_{|(\AA^n_k)^*} & \rTo^{\hat{v}_{|(\AA^n_k)^*}} &
(\AA^n_k)^* \\
\dTo^{f} & & \dTo^{\pi}\\
Bl_I=\Biproj(B) & \rTo^{v} & \PP^{n-1}_k,
\end{diagram}
where the canonical projections $f$ and $\pi$ are (trivial)
$\GG_m$-torsors 
(also called principal $\GG_m$-bundles, see \cite{MuFo82}). 
The ring $B_\pp$ is naturally $\ZZ\times\ZZ$-graded and we have the
commutative diagram
\begin{diagram}
\Proj_A(B_\pp) & \rTo & \Proj_A(B)\\
\dTo^{\hat{v}_\pp} & & \dTo^{\hat{v}}\\
\Spec(k[\Tg]_\pp) & \rTo^{\theta_\pp} & \Spec(k[\Tg])
\end{diagram}
which shows, since $\theta_\pp$ is flat, that
\begin{equation}\label{Gammap}
\Gamma(\Proj_A(B_\pp),\OO_{\Proj_A(B_\pp)})\simeq
\Gamma(\Proj_A(B),\OO_{\Proj_A(B)})_\pp. 
\end{equation}
The morphism $\Proj_A(B_\pp) \rightarrow \Proj_A(B)_{|(\AA^n_k)^*}$ is
clearly $\{1\}\times \GG_m$-equivariant, and hence, passing to the
quotient, we obtain the commutative diagram
\begin{diagram}
\Proj_A(B_\pp) & \rTo & \Proj_A(B)_{|(\AA^n_k)^*}\\
\dTo^{f_\pp} & & \dTo^{f}\\
\Biproj(B_\pp) & \rTo & \Biproj(B),
\end{diagram}
where both vertical arrows are $\GG_m$-torsors. We deduce the
isomorphism 
\begin{equation}\label{isofp}
  \OO_{\Proj_A(B_\pp)}\simeq f^*_\pp({{B_\pp}_{(0,0)}}^{\#\#}).
\end{equation}

Now we have
$\Proj_A(B_\pp)=\Spec(k[\Tg]_\pp)\times_{(\AA^n_k)^*}\Proj_A(B)_{|(\AA^n_k)^*}$
and  
$\Spec(k[\Tg]_\pp)=(\AA^n_k)^*\times_{\PP^{n-1}_k}\Spec(\OO_{\PP^{n-1}_k,s})$.
It follows that
$$\Proj_A(B_\pp)=\Spec(\OO_{\PP^{n-1}_k,s})\times_{\PP^{n-1}_k}\Proj_A(B)_{|(\AA^n_k)^*},$$
and the diagram (where both squares are commutative)
\begin{diagram}
\Proj_A(B_\pp) & \rTo & \Proj_A(B)_{|(\AA^n_k)^*}\\
\dTo^{f_\pp} & & \dTo^{f}\\
\Biproj(B_\pp) & \rTo & \Biproj(B)\\
\dTo & & \dTo^{v}\\
\Spec(\OO_{\PP^{n-1}_k,s}) & \rTo & \PP^{n-1}_k
\end{diagram}
shows that
$\Biproj(B_\pp)=Bl_I\times_{\PP^{n-1}_k}\Spec(\OO_{\PP^{n-1}_k,s})=Bl_I\times_S\Spec(\OO_{S,s})=Y$.
From this we can also deduce
$\Proj_A(B_\pp)=Bl_I\times_S\Spec((k[\Tg]/\pp)_\pp)$ where the
projection of $\Spec((k[\Tg]/\pp)_\pp)$ on $S$ is given by the
composed morphism
$\Spec((k[\Tg]/\pp)_\pp)\xrightarrow{\alpha}\Spec(\OO_{S,s})
\rightarrow S$. Moreover the localized morphism $f_\pp$ is then
$f_\pp : Bl_I\times_S \Spec((k[\Tg]/\pp)_\pp)
\xrightarrow{\mathrm{Id}\times\alpha}
Bl_I\times_S\Spec(\OO_{S,s})$. By hypothesis the scheme
$Y=Bl_I\times_S \Spec(\OO_{S,s})$ is finite on $\Spec(\OO_{S,s})$, it
is hence the spectrum of a semi-local ring and it follows (by \cite{BAC85},
proposition 5, II.5)
$${{B_\pp}_{(0,0)}}^{\#\#}\simeq \OO_{\Biproj(B_\pp)}\simeq\OO_Y.$$
By \eqref{isofp} we deduce $\OO_{\Proj_A(B_\pp)}\simeq
f^*_\pp(\OO_\Biproj(B_\pp))$. As $B_\pp$ contains an invertible
homogeneous element of non-zero degree with respect to the grading 
of $k[\Tg]$, $f_\pp$ is an affine morphism and we obtain the following
isomorphisms of $(k[\Tg]/\pp)_\pp$-modules
\begin{eqnarray*}
\Gamma(\Proj_A(B_\pp),\OO_{\Proj_A(B_\pp)}) & \simeq  &
\Gamma(\Biproj(B_\pp),{f_\pp}_*(\OO_{\Proj_A(B_\pp)})) \\
 & \simeq & \Gamma(Y,\OO_Y)\otimes_{\OO_{S,s}} (k[\Tg]/\pp)_\pp.
\end{eqnarray*}
Comparing with \eqref{Gammap}, this proves the claim \eqref{iso}.
\end{proof}

Applying this theorem with $A=k[X_1,\ldots,X_r]$ (or equivalently
$J=(0)$), we can compute the degree of the closed image of a
generically finite rational map
\begin{eqnarray*}
\PP^{r-1}_k & \xrightarrow{\lambda} & \PP^{n-1}_k\\
(X_1:\cdots:X_r) & \mapsto & (f_1:\cdots:f_n)(X_1:\cdots:X_r),
\end{eqnarray*}
if $\dim(A/(f_1,\dots,f_n))\leq 1$, and moreover decide if this 
map is generically finite by computing
$d^{r-1}-e(T,\Proj(A))$. For instance if we suppose in addition  $r=3$ 
and $n=4$ we
 recover the well-known formula to compute the degree of a surface of
 $\PP^3_k$ parameterized by 4 homogeneous polynomials of $\PP^2_k$ of
 the same degree  and without common factor 
 (see e.g. \cite{Cox01}, appendix). This theorem also yields a more general
 formula to compute the degree of the closed image of a generically
 finite rational map from an irreducible subvariety $X=\Proj(A)$ of
 $\PP^{r-1}_k=\Proj(k[X_1,\ldots,X_r])$ to $\PP^{n-1}_k$ with a finite 
 number of base points (possibly zero). For instance the closed image
 of a regular map $\lambda$ from a curve $\mathcal{C}$ of $\PP^2_k$ to
 $\PP^2_k$ given by homogeneous polynomials of same degree $d\geq 1$
 (this implies that $\lambda$ is generically finite by the second
 point of the theorem) is $d.\deg(\mathcal{C})/deg(\lambda)$.

\section{Blow-up algebras associated to a rational map}\label{BL}
Like at the beginning of the preceding section, we suppose that $k$ is a
commutative ring and $A$ is a $\ZZ$-graded 
$k$-algebra. We denote also by $\tau:k \rightarrow A_0$ the canonical
morphism of rings and consider the $k$-algebra morphism
$$\begin{array}{cccc}
h : & k[T_1,\ldots,T_n] & \longrightarrow  & A \\
& T_i & \mapsto & f_i,
\end{array}$$
where each $f_i$ is supposed to be of degree $d\geq 1$. 
 We will focus on two blow-up algebras associated to
the ideal $I=(f_1,\ldots,f_n)$ of $A$, the Rees algebra $\Rees_A(I)$ and
the symmetric algebra $\Sym_A(I)$, and show their close relation with the ideal
$\ker(h)$ of $k[T_1,\ldots,T_n]$ we would like to study. We have 
 already seen that the Rees algebra of $I$ appears naturally in
our situation, and we begin to deal with it. Hereafter, for simplicity, we will
often denote the sequence $T_1,\ldots,T_n$ by the bold letter
$\Tg$ and for instance write $k[\Tg]$ (resp.~$A[\Tg]$) instead of
$k[T_1,\ldots,T_n]$ (resp.~$A[T_1,\ldots,T_n]$).

Introducing a new indeterminate $Z$, the Rees algebra of $I$ can be
described as the image of the $A$-algebra morphism
\begin{eqnarray*}
  \beta : A[T_1,\ldots,T_n] & \longrightarrow & A[Z,Z^{-1}] \\
   T_i & \mapsto & f_iZ^{-1}.
\end{eqnarray*}
(Note that, compared with the definition \eqref{Rees} of $\beta$,  we are switching to $Z^{-1}$ in order to not overload notations in the sequel.) The kernel of $\beta$ has the following simple description that
$$\ker(\beta)=(T_1-f_1Z^{-1},\ldots,T_n-f_nZ^{-1})\cap A[\Tg].$$
This gives an easy way to compute explicitly $\ker(\beta)$ since we only need to eliminate the variable $Z^{-1}$ from the ideal $(T_1-f_1Z^{-1},\ldots,T_n-f_nZ^{-1})$. This description also shows how $\ker(\beta)$ and what is called ``moving hypersurfaces'' in different works dealing
with the implicitization problem (see \ref{overcurve} and \ref{oversurface}) are related. A moving hypersurface following the parameterization $f_1,\ldots,f_n$ is a polynomial $F\in A[\Tg]$ homogeneous in the $T_i's$ and which satisfies $F(f_1,\ldots,f_n)=0$. The latter condition implies that $F\in (T_1-f_1,\ldots,T_n-f_n)$, so moving hypersurfaces following the parameterization are exactly the homogeneous elements of $(T_1-f_1,\ldots,T_n-f_n)\cap A[\Tg]$. But it is easy to check that these homogeneous elements generate $(T_1-f_1Z^{-1},\ldots,T_n-f_nZ^{-1})\cap A[\Tg]=\ker(\beta)$ by adding the new variable $Z^{-1}$. Thus it follows that $\ker(\beta)$ is generated by the moving hypersurfaces following the parameterization.

We now turn to a second description of $\ker(\beta)$ in terms of inertia forms, which will enable us to relate it with the kernel of the map $h$.

\begin{lemma}\label{betaTF}
  $\ker(\beta)=\TF_{(Z)}((f_1-T_1Z,\ldots,f_n-T_nZ))\cap  A[T_1,\ldots,T_n]$ as ideals of the ring $A[T_1,\ldots,T_n]$.
\end{lemma}
\begin{proof}
Let us consider the $A$-algebra morphism 
$$A[T_1,\ldots,T_n,Z] \xrightarrow{\phi} A[Z,Z^{-1}] : T_i \mapsto
f_iZ^{-1},$$
and the quotient algebra $D=A[\Tg,Z]/(f_1-T_1Z,\ldots,f_n-T_nZ)$. The
kernel of the localization morphism $D \rightarrow D_Z$ is
$H^0_{(Z)}(D)$. The $A$-algebra isomorphism $D_Z\rightarrow
A[Z,Z^{-1}]$ obtained 
by sending each $T_i$ to $f_i/Z$ makes the following diagram
commutative:
\begin{diagram}
  A[\Tg,Z] & \rTo^{\phi} & A[Z,Z^{-1}] \\
  \dOnto &  & \uTo_{\wr} \\
  D & \rTo & D_Z
\end{diagram}
We hence deduce that
$\ker(\phi)=\TF_{(Z)}(f_1-T_1Z,\ldots,f_n-T_nZ)$. Since
$\ker(\beta)=\ker(\phi)\cap A[\Tg]$, the lemma follows immediately.
\end{proof}

\begin{remark}\label{TFprime} We saw that
  $H^0_{(Z)}(D)=\ker(D\rightarrow D_Z)$ and 
  is hence prime if $A$ is assumed to be an integral domain. It
  follows then easily that the ideal 
$\TF_{(Z)}(f_1-T_1Z,\ldots,f_n-T_nZ)$ is also
  prime since it is isomorphic to the kernel of the canonical map
  $A[\Tg,Z]\rightarrow D_Z$.
\end{remark}

Recall that $h$ is the map $k[T_1,\ldots,T_n]\rightarrow A$ which
sends each $T_i$ to $f_i$. If we denote by $\tau[\Tg]$ the canonical polynomial
extension 
\begin{eqnarray*}
 \tau[\Tg] : k[T_1,\ldots,T_n] & \longrightarrow & A[T_1,\ldots,T_n] \\
                        T_i & \mapsto & T_i, 
\end{eqnarray*}
which is a $k[\Tg]$-algebra morphism, we can easily check
that
\begin{eqnarray}\label{kerh}
  \ker(h) & = & \tau[\Tg]^{-1}((T_1-f_1,\ldots,T_n-f_n)) \\
   & = & \{P\in k[T_1,\ldots,T_n] : P(f_1,\ldots,f_n)=0 \} \nonumber.
\end{eqnarray} 
If we denote also by $\Theta$ the
composed morphism $i\circ\tau[\Tg]$, where $i$ is the canonical inclusion
$A[T_1,\ldots,T_n]\hookrightarrow A[T_1,\ldots,T_n,Z]$, we have the
following description of $\ker(h)$: 
\begin{lemma}\label{kerT-f}
  $\ker(h) = \Theta^{-1}(\TF_{(Z)}(f_1-T_1Z,\ldots,f_n-T_nZ))$.
\end{lemma}
\begin{proof}
  First let $P\in k[\Tg]$ such that $\Theta(P)\in
  \TF_{(Z)}(f_1-T_1Z,\ldots,f_n-T_nZ))$, that is there exists $N\in
  \NN^*$ such that $Z^N\tau[\Tg](P) \in (f_1-T_1Z,\ldots,f_n-T_nZ)$ in
  $A[\Tg,Z]$. Specializing $Z$ to $1$ we obtain $P\in\ker(h)$ by \eqref{kerh}.
  
  Now let $P\in \ker(h) \subset k[\Tg]$. As the ideal $\ker(h)$ is
  homogeneous in $k[\Tg]$, we can suppose that $P$ is
  homogeneous of degree $m$. In this case we can write  
\begin{eqnarray*}
  Z^{m}P(\Tg) & = & P(ZT_1,\ldots,ZT_n) \\
              & = & P(ZT_1,\ldots,ZT_n) - P(f_1,\ldots,f_n) \\
              & \in & (f_1-T_1Z,\ldots,f_n-T_nZ).
\end{eqnarray*}
We hence deduce that $P\in \TF_{(Z)}((f_1-T_1Z,\ldots,f_n-T_nZ))$.
\end{proof}

\begin{remark}\label{remTF} If $k\subset A_0$ then we have immediately
  $$\ker(h) = \TF_{(Z)}(f_1-T_1Z,\ldots,f_n-T_nZ)\cap k[\Tg].$$
Moreover if we suppose that $k=A_0$, $\deg(T_i)=0$ and $deg(Z)=d \geq 1$,
then we obtain $k[\Tg]=(A[\Tg,Z])_0$ and so
  $$\ker(h)=(\TF_{(Z)}(f_1-T_1Z,\ldots,f_n-T_nZ))_0.$$
\end{remark}
We are now ready to relate $\ker(h)$ and $\ker(\beta)$.
\begin{proposition}\label{h-beta} Suppose that $k\subset A_0$, then
  $$\ker(h)=\ker(\beta)\cap
  k[\Tg]=\TF_{(Z)}(f_1-T_1Z,\ldots,f_n-T_nZ))\cap k[\Tg].$$ 
Moreover, if $J$ is an ideal of $A$ such that $H^0_J(A)=0$, then we have 
$\ker(\beta)=\TF_J(\ker(\beta))$ and hence 
$$\ker(h)=\ker(\beta)\cap k[\Tg]=\TF_J(\ker(\beta))\cap k[\Tg].$$
\end{proposition}
\begin{proof}
By lemma \ref{betaTF} and remark \ref{remTF} the first claim is
proved. To prove the second part of this lemma we consider the $A$-algebra
$$D=A[\Tg,Z]/(f_1-T_1Z,\ldots,f_n-T_nZ).$$
This yields an isomorphism of $A$-algebras, that we have already used in 
lemma \ref{betaTF}, 
 $D_Z \xrightarrow{\sim}
A[Z,Z^{-1}]$ obtained by sending each $T_i$ to $f_i/Z$. Now suppose
that $J$ is an ideal of $A$ such that $H^0_J(A)=0$, then
$H^0_J(D_Z)=H^0_J(A)[Z,Z^{-1}]=0$. This implies that $H^0_J(D)
\subset H^0_{(Z)}(D)$, from we deduce  
$$  \TF_{J}((f_1-T_1Z,\ldots,f_n-T_nZ)) \subset
\TF_{(Z)}((f_1-T_1Z,\ldots,f_n-T_nZ)).$$
Now it comes easily
{\small \begin{eqnarray*}
 \TF_{(Z)}((f_1-T_1Z,\ldots,f_n-T_nZ)) & = &
 \TF_{(Z)}\TF_{J}((f_1-T_1Z,\ldots,f_n-T_nZ)) \\
 & = & \TF_{J}\TF_{(Z)}((f_1-T_1Z,\ldots,f_n-T_nZ)),
 \end{eqnarray*}}
and intersecting with $A[\Tg]$ we obtain $\ker(\beta)=\TF_J(\ker(\beta))$.
\end{proof}

This result shows the close relation between $\ker(h)$ and
$\ker(\beta)$. In fact, as we have already mentionned, a lot of information about the closed image of $\lambda$
is contained in the kernel of $\beta$ (which can be obtained by
Gr\"obner basis computations for instance). However $\ker(\beta)$ is
difficult to study in general, and we prefer 
to work with  an ``approximation'' of it which involves the study of
the symmetric algebra of $I$ that we now describe.

\medskip


Classically, we have the well known
canonical surjective morphism of $A$-algebras
\begin{eqnarray*}
  \alpha : A[T_1,\ldots,T_n]  & \longrightarrow & \Sym_A(I)
  \\
  T_i & \mapsto & f_i,
\end{eqnarray*}
whose kernel is described by
$$\ker(\alpha)=\{T_1g_1+\ldots+T_ng_n:g_i\in A[\Tg],
\sum_{i=1}^n f_ig_i=0\}.$$
Note here that $\ker(\alpha)$ is generated by  the moving hyperplanes (i.e.~moving hypersurfaces of degree 1) following the parameterization $f_1,\ldots,f_n$. 
The symmetric algebra of $I$ appears naturally by its link with
the Rees algebra of $I$ (see for instance \cite{Vas94}). We
have the following commutative diagram
\begin{diagram}
\ker(\beta) & & \lInto & & \ker(\alpha) \\ 
 & \rdTo & & \ldTo &  \\
& & A[T_1,\ldots,T_n] & & \\
 & \ldTo^{\alpha} & & \rdTo^{\beta} &  \\
\Sym_A(I) & & \rTo^{\sigma}  & & \Rees_A(I) \\ 
\end{diagram}
where $\sigma$ denotes the canonical map from $\Sym_A(I)$ to
$\Rees_A(I)$. In fact the quotient $\ker(\beta)/\ker(\alpha)$ has been
 widely studied as it gives a measure of the
difficulty in examining the Rees algebra of $I$. We
recall that the ideal $I$ is said to be of \emph{linear type} if the
canonical map $\sigma$ is an isomorphism. The following proposition
can be summarizes by 
 ``the ideal $\ker(\alpha)$ is as an approximation of the ideal $\ker(\beta)$''.
\begin{proposition}\label{reesym} Let $J$ be an ideal of $A$ such that
  the ideal  
  $I$ is of linear type outside $V(J)$ then
  $$\TF_J(\ker(\alpha))=\TF_J(\ker(\beta)).$$
If moreover $H^0_J(A)=0$ then 
$$\ker(\beta)=\TF_J(\ker(\alpha)).$$
\end{proposition}
\begin{proof}
The first assertion comes by definition: if $I$ is of linear
type outside $V(J)$ then the $A[\Tg]$-module
$\ker(\beta)/\ker(\alpha)$ is supported in $V(J)$, that is
$$J.A[\Tg] \subset \sqrt{J.A[\Tg]} \subset
\sqrt{\ann_{A[\Tg]}(\ker(\beta)/\ker(\alpha))},$$ 
which implies $\TF_J(\ker(\alpha))=\TF_J(\ker(\beta))$. The second
statement is a consequence of the first one and
proposition \ref{h-beta}. 
\end{proof}

\begin{remark}
We will see later, as a particular case, that if the first homology
group of the Koszul complex associated to the sequence 
 $f_1,\ldots,f_n$ is zero outside $V(J)$, then $I$ is of linear
type outside $V(J)$.
\end{remark}

\noindent In view of  propositions \ref{reesym} and 
\ref{h-beta}, we end this section with the following corollary which
summarizes the relations between $\ker(h)$, $\ker(\beta)$ and $\ker(\alpha)$:
\begin{corollary}\label{ker=TF} Suppose that $k\subset A_0$. If $J$ is
  an ideal of $A$ such that $H^0_J(A)=0$  
and $I$ is of linear type outside $V(J)$, 
then $$\ker(h)=\ker(\beta)\cap k[\Tg]=\TF_J(\ker(\beta))\cap
k[\Tg]=\TF_J(\ker(\alpha))\cap k[\Tg].$$ 
\end{corollary}

\section{Approximation complexes}\label{approx}
In the first part of this section we give the definition and some
basic properties 
of the \emph{approximation complexes}. These complexes was introduced in
\cite{SiVa81} and systematically developed in 
\cite{HSV82} and \cite{HSV83}.  At their most typical, they are
projective resolutions of the symmetric algebras of ideals and allow
an in-depth study of the canonical morphism $\sigma:
\Sym_A(I)\rightarrow \Rees_A(I)$. In what follows we only develop
(sometimes without proof) those properties that directly affect the
applications we are interested in. For a complete treatment on the
subject we refer the reader to the previously cited articles. In the
second part of this section we prove two new acyclicity lemmas about the
approximation complexes that we will use later to deal with the
implicitization problem. Finally, we end this section with 
 a problem introduced by David Cox in the study of the method
of ``moving surfaces'', and called \emph{Koszul
  syzygies} (see \cite{CoSc01}).


\subsection{Definition and basic properties}
Let $A$ be a commutative ring and $J$ be an ideal of $A$ generated by $r$
elements $a_1,\ldots,a_r$ (which we will often denote by the bold
letter $\ag$).
Both applications
$$ u : A[T_1,\ldots,T_r]^r \xrightarrow{(a_1,\ldots,a_r)} 
A[T_1,\ldots,T_r] : (b_1,\ldots,b_r) \mapsto \sum_{i=1}^r b_ia_i,$$
$$ v : A[T_1,\ldots,T_r]^r \xrightarrow{(T_1,\ldots,T_r)}
A[T_1,\ldots,T_r] : (b_1,\ldots,b_r) \mapsto \sum_{i=1}^r b_iT_i,$$
give two Koszul complexes $K(\ag;A[\Tg])$ and $K(\Tg;A[\Tg])$ with respective
differentials $d_\ag$ and $d_\Tg$. One can check easily that these
differentials satisfy the property 
$d_\ag\circ d_\Tg + d_\Tg\circ d_\ag =0$, and hence there exists 
three complexes, the so-called approximation complexes, which we denote
\begin{eqnarray*}
 & & \Zc_\bullet=( \ker \, d_\ag , d_\Tg) \\
 & & \Bc_\bullet=( \Im  \, d_\ag , d_\Tg) \\
 & & \Mc_\bullet=( H_\bullet(K(\ag;A[\Tg])), d_\Tg).
\end{eqnarray*}

 The $\Zc$-complex ends with the sequence $\ker(u)
\xrightarrow{v} A[T_1,\ldots,T_r] \rightarrow 0$. Since
$v(\ker(u))=\{\sum_{i=1}^r b_iT_i \ s.t. \ \sum_{i=1}^rb_ia_i=0 \}$, we deduce
that  $$H_0(\Zc)=\frac{A[T_1,\ldots,T_r]}{v(\ker(u))}\simeq\Sym_A(J).$$
A similar argument show that
$H_0(\Mc)\simeq\Sym_{A/J}(J/J^2)$. More generally one can check that
$v(\ker(u))$ annihilates the homology $A[\Tg]$-modules of $\Zc, \Bc$
and $\Mc$ which are hence modules over $\Sym_A(J)$. These homology
modules have also a nice property, which is probably one of the most
important of the approximation complexes:
\begin{proposition}\label{indepapprox}
  The homology modules of $\Zc, \Bc$ and $\Mc$ do not depend on the
  generating set chosen for the ideal $J$.
\end{proposition}
\begin{proof}
  See proposition 3.2.6 and corollary 3.2.7 of \cite{Vas94}.
\end{proof}

Some other nice properties of the approximation complexes are due to the
exactness of the standard Koszul complex $\Lc=K(\Tg;A[\Tg])$. By
definition we have the exact sequence of complexes 
\begin{equation}\label{ZLB} 0 \rightarrow \Zc \hookrightarrow \Lc
  \xrightarrow{d_\ag} 
  \Bc[-1]\rightarrow 0,
\end{equation}
where $\Bc[-1]$ denotes the translate of $\Bc$ such that
$\Bc[-1]_n=\Bc_{n-1}$. Moreover, setting $\deg(T_i)=1$ for all $i$,
the complexes $\Lc, \Zc$ and $\Bc$ (and also $\Mc$) are graded and
give for all $t\in \NN$ an exact sequence
\begin{equation*} 0 \rightarrow (\Zc)_t \hookrightarrow (\Lc)_t \xrightarrow{d_\ag}
  (\Bc[-1])_{t-1} \rightarrow 0,
\end{equation*}
where $(\Lc)_t, (\Zc)_t$ and  $(\Bc)_t$ are the part of degree $t$ of the
complexes $\Lc, \Zc$ and $\Bc$. For instance $(\Lc)_t$ is the complex of $A$-modules
$$0\rightarrow A[\Tg]_{t-r} \xrightarrow{d_\Tg} \ldots \xrightarrow{d_\Tg} {A[\Tg]}^r_{t-1}
\xrightarrow{d_\Tg} {A[\Tg]}_{t} \rightarrow 0.$$
Now from the exactness of $\Lc$ in positive degrees and the long exact
sequence associated to (\ref{ZLB}) we deduce, for all $i\geq 1$, 
an isomorphism of graded modules $H_i(\Bc) \xrightarrow{\sim}
H_i(\Zc)(1)$, where $H_i(\Zc)(1)$ denotes the graded module 
$H_i(\Zc)$ with degree shifted by 1 (i.e.~$H_i(\Zc)(1)_t=H_i(\Zc)_{t+1}$). 
 For the case $i=0$ the long exact sequence ends with
$$H_1(\Lc)=0 \rightarrow H_0(\Bc) \rightarrow  H_0(\Zc)
\xrightarrow{\pi} H_0(\Lc) \rightarrow 0,$$
hence $H_0(\Bc)\simeq \ker(\pi)$. Since $H_0(\Lc)_i=0$ if $i >
0$ and is $A$ if $i=0$, we deduce
\begin{proposition}\label{HiB}
If for all $i\geq 0$ we define the graded modules $\tilde{H}_i(\Zc)$
by $\tilde{H}_i(\Zc)_t=H_i(\Zc)_t$ for all $i\geq 1$ and all $t$,   
$\tilde{H}_0(\Zc)_t=H_0(\Zc)_t$ if $t\neq 0$ and
$\tilde{H}_0(\Zc)_0=0$, then we have an isomorphism of graded modules
$$H_i(\Bc) \xrightarrow{\sim} \tilde{H}_i(\Zc)(1) \ \mathrm{for \
  all} \ i\geq 0.$$
\end{proposition}

We have another natural exact sequence of complexes involving the complex
$\Mc$:
\begin{equation}\label{BZM}
  0\rightarrow \Bc \rightarrow \Zc \rightarrow \Mc \rightarrow 0.
\end{equation}
Proposition \ref{HiB} and the long exact sequence associated to
(\ref{BZM}) give the following graded long exact sequence:
\begin{eqnarray}\label{lZMB}
 \ldots \rightarrow \tilde{H}_i(\Zc)(1) \rightarrow H_i(\Zc)
\rightarrow H_i(\Mc) \rightarrow \tilde{H}_{i-1}(\Zc)(1) \rightarrow
\ldots \\
 \ldots \rightarrow \tilde{H}_0(\Zc)(1) \rightarrow H_0(\Zc)
\rightarrow H_0(\Mc) \rightarrow 0. \nonumber
\end{eqnarray}
This sequence yields two results on these approximation
complexes. The first one is contained in the following proposition.
\begin{proposition}\label{MtoZ} Suppose that $A$ is a noetherian ring and $i\geq 1$. If
  $H_i(\Mc)=0$ then  $H_i(\Zc)=0$.  In particular, if $\Mc$ is acyclic
  then $\Zc$ is also acyclic.
\end{proposition}
\begin{proof} From the long exact sequence (\ref{lZMB}) we have, for
  all $i\geq 1$, the exact sequence
  $$\tilde{H}_i(\Zc)(1) \rightarrow H_i(\Zc) \rightarrow  H_i(\Mc).$$
  By hypothesis $H_i(\Mc)=0$, and since $i\geq 1$, we deduce a
  surjective morphism of $A[\Tg]$-modules $H_i(\Zc)(1) \rightarrow
  H_i(\Zc)$. Now $A$ is noetherian, so $A[\Tg]$ is also and we deduce
  that $H_i(\Zc)$ is a $A[\Tg]$-module of finite type. This implies
  that our surjective morphism is in fact bijective. Taking into
  account the degrees we obtain isomorphisms $ H_i(\Zc)_{t+1}
  \xrightarrow{\sim} 
  H_i(\Zc)_t$ for all $t$. As $H_i(\Zc)_{-1}=0$ we deduce by iteration
  that $H_i(\Zc)_{t}=0$ for all $t$.
\end{proof}

\begin{remark}\label{remMtoZ} Notice that, by definition of $\Mc$, if $H_i(K(\ag;A[\Tg]))=0$
   then $H_i(\Mc)=0$.
\end{remark}

The second result given by (\ref{lZMB}) is a condition so that the
canonical morphism $\sigma:\Sym_A(J)\rightarrow \Rees_A(J)$ becomes an
isomorphism: in this case we say that $J$ is of \emph{linear type}. By
definition we have 
$$\tilde{H}_0(\Zc)(1)=\oplus_{t\geq 1} H_0(\Zc)_t = \oplus_{t\geq 1}
\Sym_A^t(J) =\Sym_A^+(J),$$
and hence, examining the end of (\ref{lZMB}), we obtain the exact
sequence
$$ H_1(\Mc)\rightarrow \Sym_A^+(J) \xrightarrow{\mu} \Sym_A(J) \rightarrow  
\Sym_{A/J}(J/J^2) \rightarrow 0.$$
The map $\mu$ is called the \emph{downgrading} map (see \cite{Vas94},
chap. 3), it is induced by the morphism
$$ T_A^n(J) \rightarrow T_A^{n-1}(J) : x_1\otimes\ldots\otimes x_n
\mapsto x_1x_2\otimes x_3\otimes\ldots\otimes x_n,$$
where $T_A^n(J)$ denotes the degree $n$ part of the tensor algebra
associated to $J$ over $A$. Denoting  $\gr_J(A):=\oplus_{i\geq 0}
J^i/J^{i+1}$ and  $\Rees_A^+(J):=\oplus_{i\geq 1} J^i$, we obtain the
 commutative diagram 
\begin{diagram}
 H_1(\Mc) & \rTo & \Sym_A^+(J) & \rTo^{\mu} &  \Sym_A(J)
 & \rTo^{\pi} & \Sym_{A/J}(J/J^2) & \rTo & 0 \\
   & & \dTo & & \dTo_{\sigma} & & \dTo_{\gamma} & & \\
 0 & \rTo & \Rees_A^+(J) & \rTo &  \Rees_A(J)
 & \rTo & \gr_{J}(A) & \rTo & 0, \\
\end{diagram}
which gives the following proposition. 
\begin{proposition}\label{lintype}
 If $H_1(\Mc)=0$ then $ \Sym_A(J)\simeq \Rees_A(J)$, that is $J$ is of
 linear type.
\end{proposition}
\begin{proof}
 For all $i\geq 0$ we have the commutative diagram
\begin{diagram}
0 & \rTo & \Sym_A^{i+1}(J) & \rTo^{\mu} &  \Sym_A^i(J) \\
   & & \dTo_{\sigma_{i+1}} & & \dTo_{\sigma_i} \\
 0 & \rTo & J^{i+1} & \rTo &  J^i. \\   
\end{diagram}
But $\sigma_0$ is an isomorphism, which implies that $\sigma_1$ is
injective. As $\sigma_1$ is already surjective it is an isomorphism. By
iteration we show that $\sigma_t$ is an isomorphism for all $t\in
\NN$, and hence that $
\Sym_A(J)\simeq\Rees_A(J)$.
 \end{proof}

 This result gives a criterion for an ideal to be of \emph{linear
   type}, that is $\sigma$ is an isomorphism. If the ring $A$ is
 noetherian one can show that $\sigma$ is an isomorphism if and only if
 $\gamma$ is (see \cite{Vas94}, theorem 2.2.1). There also exists
 other criterions for an ideal
 to be of linear type based on some properties of a set of generators
 of $J$ (see \cite{Vas94}, paragraph 3.3) but we will not go further
 in this direction here.

 
\subsection{Acyclicity criteria}
We will here prove two acyclicity
results that we will apply later to the implicitization
problem. We first recall the acyclicity criterion of Peskine-Szpiro:
\begin{lemma}\label{acylem} {\rm (Acyclicity lemma)}
Suppose that $R$ is a commutative noetherian ring, $I$ a proper ideal, and
$$C_\bullet : 0\rightarrow M_n \rightarrow M_{n-1}\rightarrow \ldots
\rightarrow M_0$$
a complex of finitely generated $R$-modules with homology
$H_k=H_k(C_\bullet)$ such that each non-zero module $H_k$, for all
$k=1,\ldots,n$, satisfies $\depth_I(H_k)=0$. Then $\depth_I(M_k)\geq
k$ for all $k=1,\ldots,n$ implies  $C_\bullet$ is exact.
\end{lemma}
\begin{proof} We refer to \cite{BuEi73}, lemma 3.
\end{proof}

 Hereafter we assume that $A$ is a noetherian commutative
 $\NN$-graded ring, and we denote by $\mm$
its ideal generated by elements of strictly positive degree, that is
$\mm=A_+=\oplus_{\nu>0}A_\nu$. The following proposition is our first
acyclicity result. 
\begin{proposition}\label{acycwithoutbp} Let $I=(a_1,\ldots,a_n)$ be
  an ideal of $A$ such that both ideals $I$ and $\mm$ have the same
  radical, and suppose that  $\sigma=\depth_\mm(A)\geq 1$. Then the
  homology modules $H_i(\Zc)$ vanish for $i\geq \max(1,n-\sigma)$,
  where $\Zc$ denotes the $\Zc$-complex associated to the ideal $I$.\\
  In particular if $n\geq 2 $ and $\sigma=\depth_\mm(A) \geq n-1$ then
  the complex $\Zc$ is acyclic. 
\end{proposition}
\begin{proof} This proof is based on chasing depths and the acyclicity
  lemma \ref{acylem}. Let $K_\bullet(\ag;A)$ be the Koszul complex
  associated to the sequence $a_1,\ldots,a_n$ in $A$. We denote by
  $H_i(\ag;A)$ its homology modules and by $Z_i$
  (resp.~$B_i$)  the $i$-cycles (resp.~the $i$-boundaries) of
  $K_\bullet(\ag;A)$. The $\Zc$-complex is of the form
   $$ 0 \rightarrow \Zc_{n-1} \rightarrow \Zc_{n-2} \rightarrow \ldots
   \rightarrow \Zc_1 \rightarrow \Zc_0 \rightarrow 0.$$
   Since $I$ and $\mm$ have the same radical the homology modules
   $H_i(\ag;A)$ are supported on $V(\mm)$ for $I$ annihilates them. Hence, by proposition \ref{MtoZ} and remark \ref{remMtoZ}, homology modules 
   $H_i(\Zc)$ are 
   supported on $V(\mm)$ for all $i\geq 1$, and
   consequently we have 
   \begin{equation}\label{p2}
     \depth_\mm(H_i(\Zc))=0 \ \ {\rm for \ all} \  H_i(\Zc)\neq 0 , \ i\geq 1.  
   \end{equation}

   Always since $I$ and $\mm$ have the same radical we
  have $\sigma=\depth_\mm(A)=\depth_I(A)$ and hence deduce by lemma \ref{acylem} that 
  $H_i(\ag;A)=0$ for $i>n-\sigma$. It follows that the following
  truncated Koszul complex is exact:
  $$0\rightarrow K_n \rightarrow K_{n-1} \rightarrow \ldots \rightarrow
  K_{n-\sigma+1} \rightarrow K_{n-\sigma},$$
  and, since $\sigma\geq 1$, we have exact sequences
  \begin{eqnarray*}
     &  & 0 \rightarrow K_n \rightarrow B_{n-1} \rightarrow 0, \\
    &  & 0 \rightarrow B_{n-1} \rightarrow K_{n-1} \rightarrow B_{n-2}
    \rightarrow 0,\\ 
      & & \hspace{1cm} \vdots \hspace{2.7cm} \vdots \\
     & &  0 \rightarrow B_{n-\sigma+1} \rightarrow K_{n-\sigma+1} \rightarrow
     B_{n-\sigma} \rightarrow 0.
   \end{eqnarray*}
We can now use standard properties of depth (see \cite{Eis94} corollary 18.6). 
   As $\depth_\mm(K_i)\geq \sigma$ for all $i$, it follows by iterations
   that $\depth_\mm(B_i) \geq \sigma-(n-i)+1$ for $i$  from
   $n-\sigma$ to $n-1$. Since $\sigma\geq 1$ we have $Z_n=0$ and
   $Z_i=B_i$ for $n-1 \geq i > n-\sigma$. Moreover $Z_{n-\sigma}
   \subset K_{n-\sigma}$ and  $\depth_\mm(K_{n-\sigma})\geq 1$ hence
   $\depth_\mm(Z_{n-\sigma})\geq 1$. We finally obtain, using the fact that $\Zc_i=Z_i[\Tg]$ for all $i=0,\ldots,n$,
   \begin{equation}\label{p1}
     \Zc_n=0 \ \ {\rm and} \ \ \depth_\mm(\Zc_i)\geq i-(n-\sigma-1) \
     \ {\rm for} \ \ n > i \geq n-\sigma. 
   \end{equation}
   By lemma \ref{acylem}, (\ref{p2}) and (\ref{p1}) show that the
   complex
$$ 0 \rightarrow \Zc_{n-1} \rightarrow \ldots \rightarrow
\Zc_{n-\sigma+1} \rightarrow \Zc_{n-\sigma-1}$$ 
   is exact, i.e.~$H_i(\Zc)=0$ for $i\geq \max(1,n-\sigma)$.
\end{proof}

The preceding result is valid only if both ideals $I$ and $\mm$ have the
same radical. Consequently, when we will apply this result to the
implicitization problem, it will be useful only in the absence of
base points. In order to deal with certain cases where there exists base
points we prove a second acyclicity lemma. 

For any ideal $J$ of a ring $R$ we denote by $\mu(J)$ the minimal
number of generators of $J$. We recall the standard definition of a
projective local complete intersection ideal.
\begin{definition}\label{lcidef}
Let $I$ be an ideal of $A$. $I$ is said to be a local complete
intersection in $\Proj(A)$ if and only if  
for all $\mathfrak{p}\in \Spec(A)\setminus V(\mm)$ we have $\mu(I_\mathfrak{p})=\depth_{I_\mathfrak{p}}(A_\mathfrak{p})$.
\end{definition}

\begin{proposition}\label{lciprop} Let $I=(a_1,\ldots,a_n)$ be an
  ideal of $A$ such  that $I$ is a 
  local complete intersection in $\Proj(A)$ and $n\geq 2$. Suppose  
  $\depth_\mm(A)\geq n-1$ and $\depth_I(A)=n-2$, then the
  $\Zc$-complex associated to $I$ is acyclic. 
\end{proposition}
\begin{proof} Notice first that, since $\depth_I(A)=n-2$, the minimal number of generators of $I$ is at least $n-2$. If $\mu(I)=n-2$ then by
   proposition \ref{indepapprox} we obtain immediately that $\Mc$ is
   acyclic and so $\Zc$, since $I$ is a complete intersection
   ideal. Consequently we suppose that $\mu(I)\geq n-1$. 
   The proof is now very similar to the one of proposition
  \ref{acycwithoutbp} and is also based on chasing depths and the acyclicity
  lemma \ref{acylem}. We keep the notation of proposition
  \ref{acycwithoutbp}, that is $K_\bullet(\ag;A)$ denotes the Koszul complex
  associated to the sequence $a_1,\ldots,a_n$ in $A$, for all $i\geq
  0$   $H_i(\ag;A)$ denote its homology modules and $Z_i$
  (resp.~$B_i$)  its $i$-cycles (resp.~its $i$-boundaries). The
  $\Zc$-complex is of the form 
   $$ 0 \rightarrow \Zc_{n-1} \rightarrow \Zc_{n-2} \rightarrow \ldots
   \rightarrow \Zc_1 \rightarrow \Zc_0 \rightarrow 0.$$
   Since $I$ is a local complete intersection in $\Proj(A)$, proposition
   \ref{indepapprox} shows   
   that the homology modules $H_i(\Mc)$ for all $i\geq 1$ are
   supported on $V(\mm)$, 
   where $\Mc$ is the approximation complex associated to the homology
   of the ideal $I$. Now by proposition \ref{MtoZ} it follows that the
   homology modules $H_i(\Zc)$ are  also 
   supported on $V(\mm)$ for all $i\geq 1$, and
   consequently we have 
   \begin{equation}\label{plci2}
     \depth_\mm(H_i(\Zc))=0 \ \ {\rm for \ all } \ H_i(\Zc)\neq 0, \ i\geq 1.  
   \end{equation}

   Since we have supposed $\depth_I(A)=n-2$, the
   homology modules  $H_i(\ag;A)$ vanish for $i>2$. The
    following truncated Koszul complex is hence exact:
  $$0\rightarrow K_n \rightarrow K_{n-1} \rightarrow \ldots \rightarrow
  K_{3} \rightarrow K_{2},$$
  and we have exact sequences
  \begin{eqnarray*}
     &  & 0 \rightarrow K_n \rightarrow B_{n-1} \rightarrow 0, \\
    &  & 0 \rightarrow B_{n-1} \rightarrow K_{n-1} \rightarrow B_{n-2}
    \rightarrow 0,\\ 
      & & \hspace{1cm} \vdots \hspace{2.7cm} \vdots \\
     & &  0 \rightarrow B_{3} \rightarrow K_{3} \rightarrow
     B_{2} \rightarrow 0.\\
   \end{eqnarray*}
   As $\depth_\mm(K_i)\geq n-1$ for all $i=1,\ldots,n$, it follows by
   iterations 
   that $\depth_\mm(B_i) \geq i$ for $i$  from
   $2$ to $n-1$. But we have $Z_n=0$ and
   $Z_i=B_i$ for $n-1 \geq i > 2$, so $\depth_\mm(Z_i)\geq i$ for
   $n-1 \geq i\geq 3$. We have also $Z_1\subset K_1$ and
   $\depth(K_1)\geq 1$, so we deduce $\depth_\mm(Z_1)\geq 1$.\\
   Now since $\Zc_i=Z_i[\Tg]$ we obtain
   \begin{equation}\label{plci1}
     \Zc_n=0 \ \ {\rm and} \ \ \depth_\mm(\Zc_i)\geq i \
     \ {\rm for} \ \ i \geq 3 \ \ {\rm and} \ i=1. 
   \end{equation}
   If we prove that $\depth(\Zc_2)\geq 2$, (\ref{plci2}) and
   (\ref{plci1}) imply that the 
   $\Zc$-complex associated to $I$ is exact by lemma \ref{acylem}.
   To show this last inequality we prove that in both cases
   $\mu(I)=n-1$ 
   and $\mu(I)=n$ we have $H_2(\Mc)=0$. For this we need the two
   lemmas stated just 
   after this proposition. If $\mu(I)=n-1$, lemma \ref{lcilemma}
   implies immediately 
   $H_2(\Mc)=0$. If $\mu(I)= n$  lemma \ref{lcilemma} implies that
   $H_2(\Mc)=H^0_\mm(H_2(\Mc))$, and by lemma
   \ref{tech} we obtain that $H^0_\mm(H_2(\Mc))=0$, so $H_2(\Mc)=0$. 
   It follows that in both cases $\Bc_2=\Zc_2$, and since we have proved
   $\depth_\mm(B_2)\geq 2$, we deduce  $\depth_\mm(\Bc_2)\geq 2$, and
   finally $\depth_\mm(\Zc_2)\geq 2$.
    \end{proof}

 \begin{lemma}\label{tech} Let $I=(a_1,\ldots,a_n)$ be an ideal of $A$
  such that $\depth_\mm(A) > \depth_I(A)$, then
  $H^0_\mm(H_{n-\depth_I(A)}(\ag;A))=0$. 
\end{lemma}
\begin{proof} This lemma is just a standard use of the classical spectral
  sequences associated to the double complex 
$$\begin{array}{ccccccccccc}
   0 & \rightarrow & \CC^0_\mm (\wedge^nA^n) &  \xrightarrow{d_\ag} & \ldots &
   \xrightarrow{d_\ag} & \CC^0_\mm (\wedge^1A^n) & \xrightarrow{d_\ag} & \CC^0_\mm(A) &
   \rightarrow & 0 \\  
    & & \downarrow &  &   &  & \downarrow &  & \downarrow &  &  \\
  0 & \rightarrow & \CC^1_\mm (\wedge^nA^n) &  \xrightarrow{} & \ldots
  & \xrightarrow{} & 
  \CC^1_\mm (\wedge^1A^n) & \xrightarrow{} & \CC^1_\mm (A) & \rightarrow & 0 \\
    & & \downarrow &  &  &  & \downarrow &  & \downarrow &
    &  \\
    & & \vdots &  & \vdots &  & \vdots &  & \vdots &
    &  \\
     & & \downarrow &  &  &  & \downarrow &  & \downarrow &
    &  \\
    0 & \rightarrow & \CC^{r}_\mm (\wedge^nA^n) &  \xrightarrow{} & \ldots & \xrightarrow{} &
  \CC^{r}_\mm (\wedge^1A^n) & \xrightarrow{} & \CC^{r}_\mm (A) &  \rightarrow & 0 .\\ 
\end{array}$$
Its first row is the Koszul complex associated to the sequence
$a_1,\ldots,a_n$,  and its columns are the classical \v{C}ech complexes
(recall that $C^0_\mm(A)=A$, and if $x_1,\ldots,x_r$ is a system of
parameters of $A$, then for all 
$t\geq 1$ we  have $\CC^t_\mm(A)=\bigoplus_{1\leq i_1 < \cdots <  i_t \leq r} A_{x_{i_1}x_{i_2}\cdots x_{i_t}}$).  We
know that $H_i(\ag;A)=0$ for $i>n-\depth_I(A)$, and also that
$H^i_\mm(A)=0$ if $i < \depth_\mm(A)$ which implies 
$H^i_\mm(A)=0$ for all $i \leq \depth_I(A)$. Examining the two
filtrations by rows and by columns we deduce that
$H^0_\mm(H_{n-\depth_I(A)}(\ag;A))=0$. 
\end{proof}

\begin{lemma}\label{lcilemma} Let $I=(a_1,\ldots,a_n)$ be an ideal of
  $A$ and let $\Mc$ be the $\Mc$-complex associated to $I$. We have
  the two following properties:
  \begin{itemize}
    \item[a)]  For all $i > \zeta:=\mu(I) -\depth_I(A)$ we have $H_i(\Mc)=0$.
     \item[b)] If we suppose that for all $\mathfrak{p}\in
Spec(A)\setminus V(\mm)$ we
have $\zeta > \zeta_\mathfrak{p}:=\mu(I_\mathfrak{p})
-\depth_{I_\mathfrak{p}}(A_\mathfrak{p})$, then 
$H_\zeta(\Mc)=H^0_\mm(H_\zeta(\Mc))$.
\end{itemize}
\end{lemma} 
\begin{proof}
This lemma is a direct consequence of the proposition
\ref{indepapprox}. Indeed, for the first statement, we can construct
the complex $\Mc$ from the  
Koszul complex of a sequence of elements $a'_1,\ldots,a'_{\mu(I)}$
which generate $I$. This Koszul complex has all its homology groups 
$H_p(\ag';A)=0$ for $p>\mu(I)-\depth_I(A)$, and we deduce, by
definition of $\Mc$, that $H_p(\Mc)=0$ for such $p$. Now for all
$\mathfrak{p}\in \Spec(A)$ and all $i\geq 0$ we have
$H_i(\Mc(I))_\mathfrak{p}\simeq H_i(\Mc(I)_\mathfrak{p})\simeq
H_i(\Mc(I_\mathfrak{p}))$, where the last isomorphism is true by 
proposition \ref{indepapprox}. From the first statement of this lemma
we deduce that $H_\zeta(\Mc(I_\mathfrak{p}))=0$ for all
$\mathfrak{p}\notin V(\mm)$, since $\zeta > \zeta_\mathfrak{p}$. It
follows that $H_\zeta(\Mc)_\mathfrak{p}=0$ for all
$\mathfrak{p}\notin V(\mm)$, which implies
$H_\zeta(\Mc)=H^0_\mm(H_\zeta(\Mc))$. 
\end{proof}

\subsection{Syzygetic ideals and Koszul syzygies}
We will  now investigate a problem introduced by David Cox, always in
relation with the implicitization problem, 
and solved by himself and Hal Schenck in
\cite{CoSc01}, theorem 1.7, in $\PP^2$. Being given a homogeneous ideal
$I=(a_1,a_2,a_3)$ of 
$\PP^2$ of codimension two, the authors showed that the module of
syzygies of $I$ vanishing at the scheme locus 
$V(I)$ is generated by the Koszul sygyzies if and only if $V(I)$ is a
local complete intersection in $\PP^2$. In other words, if $Z_1$
(resp.~$B_1$) denotes the first syzygies (resp.~the first boundaries)
of the Koszul complex associated to the sequence $a_1,a_2,a_3$, they
showed that $Z_1\cap \TF_\mm(I).A^n = B_1$ if and only if $I$ is a local
complete intersection in $\PP^2$, where $\mm$ is the irrelevant ideal
of the ring $A=k[X_1,X_2,X_3]$ defining $\PP^2$. This result is
interesting in at least two points: the first one is that it gives an
algorithmic way to test if the ideal $I$ is a local complete
intersection of $\PP^2$ or not, and the second is that it is the key point of
 a new formula for the implicitization problem (see e.g. 
\cite{BCD02}). 

In what follows we show that this problem is closely related to
syzygetic ideals and generalize the preceding result. Let $A$ be
a $\NN$-graded commutative noetherian ring, $\mm=\oplus_{\nu>0}A_\nu$ be its
irrelevant ideal,  
and let $I$ be another ideal of $A$ generated by $n$ 
elements, say  $a_1,\ldots,a_n$. We denote by $Z_1$
(resp.~$B_1$) the first syzygies (resp.~the first boundaries)
of the Koszul complex over $A$ associated to sequence
$a_1,\ldots,a_n$. Our aim is to formulate a necessary and sufficient
condition so that the inclusion 
$$ B_1 \subset Z_1\cap \TF_\mm(I).A^n$$
becomes an equality (notice that we always have $ B_1 \subset Z_1\cap
I.A^n$ and also $Z_1\cap I.A^n \subset Z_1\cap \TF_\mm(I).A^n$).

    We have an exact sequence of $A$-modules:
$$0 \rightarrow B_1 \hookrightarrow I.A^n \rightarrow \Sym_A^2(I)
\rightarrow 0,$$
where the last map sends $(x_1,\ldots,x_n) \in I.A^n$ to the equivalence
class of $(x_1\otimes a_1+\ldots + x_n\otimes a_n)$ in  
$\Sym_A^2(I)$. We have the other exact sequence
$$0\rightarrow Z_1\cap I.A^n \rightarrow I.A^n \rightarrow I^2
\rightarrow 0,$$
where the last map sends $(x_1,\ldots,x_n) \in I.A^n$ to $x_1a_1+\ldots
+ x_na_n \in I^2$. Putting these two exact sequences together we
obtain the following commutative diagram
\begin{diagram}
0 & \rTo & B_1 & \rTo & I.A^n & \rTo & \Sym_A^2(I) & \rTo & 0 \\
  &      & \dInto     &      & \dTo_{\wr} &             & \dOnto & \\
0 & \rTo & Z_1\cap I.A^n & \rTo & I.A^n & \rTo & I^2 & \rTo & 0, \\
\end{diagram}
where the vertical maps are the canonical ones. From the snake lemma
we deduce immediately the isomorphism
$$\frac{Z_1\cap I.A^n}{B_1}\simeq \ker(\Sym_A^2(I)
\rightarrow I^2).$$
Following \cite{SiVa81}, we denote $\ker(\Sym_A^2(I)
\rightarrow I^2)$ by $\delta(I)$ (even if it is not the real
definition of this invariant of $I$ but a property) and will say that $I$ is
\emph{syzygetic} if $\delta(I)=0$. We just proved:
\begin{lemma}
$B_1=Z_1\cap I.A^n$ if and only if $I$ is syzygetic.
\end{lemma}
We state now the main result of this paragraph.
\begin{theorem}\label{theokoszul} If $\depth_\mm(A) >
  \depth_I(A)=n-1$ then the following statements are equivalent:
\begin{itemize}
\item[a)] $B_1=Z_1\cap I.A^n$
\item[b)] $B_1=Z_1\cap \TF_\mm(I).A^{n}$.
\item[c)] $I$ is syzygetic.
\end{itemize}
\end{theorem}
\begin{proof}
We have already prove that a) is equivalent to c), and 
as $$B_1 \subset Z_1\cap I.A^n \subset Z_1\cap
\TF_\mm(I).A^{n},$$ it remains only to show that a) implies b).

Let $(x_1,\ldots,x_n) \in Z_1\cap\TF_\mm(I).A^n$. By
definition there exists an integer $\nu$ such that
$$\mm^\nu(x_1,\ldots,x_n) \subset Z_1\cap I.A^n=B_1.$$
This last implies that for all $\xi \in \mm$ the equivalence class of
$\xi^\nu.(x_1,\ldots,x_n)$ in $H_1(\ag,A)$ is 0. But by lemma
\ref{tech} we have $H^0_\mm(H_1(\ag,A))=0$, and we deduce that the
equivalence class of $(x_1,\ldots,x_n)$ in $H_1(\ag,A)$ is 0, that is 
$(x_1,\ldots,x_n) \in B_1$. 
\end{proof}

We now give a more explicit condition on the ideal $I$ to obtain
the desired property $B_1=Z_1\cap \TF_\mm(I).A^n$. For this we 
consider the approximation 
complex $\Mc$  attached to the ideal
$I$. From now we suppose that $A=k[X_1,\ldots,X_n]$ with $k$ a
noetherian commutative ring, 
$I=(a_1,\ldots,a_n)$, and we take 
$\mm=(X_1,\ldots,X_n)$ (in this case we have
$\depth_\mm(A)=n$). 

\begin{proposition} Suppose $k$ is a Cohen-Macaulay ring.
  If $I$ is a local complete intersection in   
  $\PP^{n-1}_k$ and $\codim(I)=n-1$, then
  $B_1=Z_1\cap\TF_\mm(I).A^n$. 
\end{proposition}
\begin{proof} Since $\depth_I(A)=n-1$, we have only to prove that $I$
  is syzygetic by theorem \ref{theokoszul}. Now lemma \ref{lcilemma}
  shows that $H_1(\Mc)=0$ in both cases $\mu(I)=n-1$ and $\mu(I)=n$:
  if $\mu(I)=n-1$ we have directly $H_1(\Mc)=0$ by a) and, if
  $\mu(I)=n$,  b) implies $H_1(\Mc)=H^0_\mm(H_1(\Mc))$ and lemma \ref{tech} shows
  that this last is zero. By proposition  \ref{lintype} it follows
  that $I$ is of linear type, which implies that $I$ is syzygetic. 
\end{proof}

In the case $n=3$ the preceding proposition becomes an equivalence,
due to a particular property of $\delta(I)$ which is stated in theorem
2.2 of \cite{SiVa81}. The following proposition generalizes the 
theorem 1.7 of \cite{CoSc01}, where $k$ is a field, to the case 
$k$ is a noetherian regular ring.   
\begin{proposition}\label{coxschenck} Suppose $k$ is a 
  regular ring, $n=3$ and $\codim(I)=2$. Then $I$ is
  a local complete 
  intersection in $\PP^2_k$ if and only if
  $B_1=Z_1\cap\TF_\mm(I).A^3$. 
\end{proposition}
\begin{proof} The only point we have to prove is that, under our
  hypothesis, if $I$ is syzygetic then $I$ is a local complete
  intersection in $\PP^2_k$. For this we will use theorem 2.2 of
  \cite{SiVa81}: let $R$ be a noetherian local ring of depth two and let $P$ an
  ideal of codimension two and projective dimension one, then $P$ is
  generated by a regular sequence if and only if $\delta(P)=0$.

  Clearly $I_\mathfrak{p}\simeq \TF_\mm(I)_\mathfrak{p}$ for all $\mathfrak{p}
  \notin V(\mm)$. As the invariant $\delta(I)$ of $I$ 
  satisfies $\delta(I)_\mathfrak{p} \simeq \delta(I_\mathfrak{p})$ for
  all $\mathfrak{p} \in \Spec(A)$, if we prove that $TF_\mm(I)$ has
  projective dimension 1 then theorem 2.2 of \cite{SiVa81} implies that
  $\TF_\mm(I)$ is a local complete intersection in $\PP^2_k$, which is
  equivalent to say that $I$ is.

  From the definition of $\TF_\mm(I)$, it follows that $\mm$ is not
  associated to $\TF_\mm(I)$ (it is easy to see that $\TF_\mm(I) : \mm =
  \TF_\mm(I)$) and hence $$\depth_\mm(A/\TF_\mm(I)) > 0.$$
  But we know that $\depth_\mm(A/\TF_\mm(I)) \leq
  \dim(A/\TF_\mm(I))$, and we know also that $\dim(A/\TF_\mm(I))=\dim(A) -
  \codim(\TF_\mm(I))=1$. 
  We deduce $$\dim(A/\TF_\mm(I)) - \depth_\mm(A/\TF_\mm(I))=0.$$
  Now, since $k$ is assumed to be regular, the Auslander-Buchsbaum formula gives 
  $$\pd(A/\TF_\mm(I))=\depth_\mm(A) - \depth_\mm(A/\TF_\mm(I)).$$
  As $\depth_\mm(A)=\dim(A)=\dim(A/\TF_\mm(I))+\codim(\TF_\mm(I))$ we
  deduce  $\pd(A/\TF_\mm(I))=\codim(\TF_\mm(I))=2$, 
  whence $\pd(\TF_\mm(I))=1$. 
\end{proof}

\section{The implicitization problem via the approximation complexes}

As we have seen, the kernel of the map $h:k[\Tg] \rightarrow A$ which sends $T_i$ to $f_i\in
A_d$ for $i=1,\ldots,n$, defines the closed image of its associated
map $\lambda$ \eqref{deflambda}. In this section we will investigate the
so-called implicitization problem. Hereafter we suppose that $k$ is a
field and that $A$ is the polynomial ring $k[X_1,\ldots,X_{n-1}]$,
with $n\geq 3$, whose 
irrelevant ideal is denoted $\mm=(X_1,\ldots,X_{n-1})\subset A$. If
the morphism $\lambda$ is generically 
finite then the closed image of $\lambda$ is a hypersurface  
of $\PP^{n-1}_k$, and the implicitization problem consists in
computing explicitly its equation (up to a nonzero multiple in $k$). In what
follows we will apply techniques using the approximation complexes we
have introduced in section \ref{approx}. We will denote by $I$ the
ideal of $A$ generated by the polynomials $f_1,\ldots,f_n$ which are
supposed to be of the same degree $d\geq 1$, considering each $X_i$ of
degree 1 in $A$. We will also denote by $\Zc$ and $\Mc$ the two
approximation complexes associated to $I$. The basic idea is to show
that the $\Zc$-complex is acyclic under certain conditions, 
giving thus a resolution of the symmetric algebra $\Sym_A(I)$, and then
obtain the implicit equation as the determinant of some of its homogeneous
components. Before going further into different particular
cases we give the main ingredient of this section.  

The symmetric algebra of $I$ is naturally
bi-graded by the exact sequence:
$$0\rightarrow \ker(\alpha) \rightarrow A[T_1,\ldots,T_n] \xrightarrow{\alpha}
\Sym_A(I) \rightarrow 0,$$
where $\ker(\alpha)=(\sum b_iT_i \ | \ b_i \in A \ \sum b_if_i=0
)\subset A[T_1,\ldots,T_n]$. We denote by $\Sym_A(I)_\nu$ the
graded part of $\Sym_A(I)$ corresponding to the grading of $A$,
that is, to be more precise, $\Sym_A(I)_\nu=\oplus_{l\geq 0} A_\nu\Sym_A^l(I)$.

\begin{proposition}\label{annprop} Suppose that $I$ is of
  linear type outside $V(\mm)$,  and let $\eta$ be an integer 
  such that $H^0_\mm(\Sym_A(I))_\nu=0$  for 
  all $\nu\geq \eta$ then 
  $$\ann_{k[\Tg]}(\Sym_A(I)_\nu)=\ker(h) \ \ {\rm for \ all}  \ \nu \geq
  \eta.$$
\end{proposition}
\begin{proof} By definition, for $\nu \in \NN$,
  $$\ann_{k[\Tg]}(\Sym_A(I)_\nu)=
  \{f\in k[\Tg] \, : \, f.\Sym_A(I)_\nu=0 \},$$ and from the description
  of $\Sym_A(I)$ we deduce
  \begin{equation}\label{hypann}
  \ann_{k[\Tg]}(\Sym_A(I)_\nu)=
  \{f\in k[\Tg] \, : \, f.A_\nu[\Tg]\subset \ker(\alpha) \}.
   \end{equation}
  In the same way, $H^0_\mm(\Sym_A(I))_\nu=0$ is equivalent to 
  \begin{equation}\label{hypH}
 \{ f\in A_\nu[\Tg] \, : \, \exists n \, \mm^nf\subset \ker(\alpha)
 \}={\ker(\alpha)}_\nu.
  \end{equation}
  From (\ref{hypann}) and (\ref{hypH}) we easily check that
  {\small \begin{eqnarray}\label{ann}
    \ann_{k[\Tg]}(\Sym_A(I)_\nu)& = &\{f\in k[\Tg] \, : \,
     f.A_\eta[\Tg]\subset \ker(\alpha) \}  \\ \nonumber
     & = & \{f\in k[\Tg] \, : \,
    f.A_{\geq \eta}[\Tg]\subset \ker(\alpha) \}, \nonumber
  \end{eqnarray}}
for all $\nu \geq  \eta$.

Now since $I$ is of linear type outside $V(\mm)$ and
$H^0_\mm(A)=0$, corollary \ref{ker=TF}
   implies that $\ker(h)$ is exactly $\TF_\mm(\ker(\alpha))\cap
   k[\Tg]$ (observe here that $\TF_\mm=\TF_{\mm[\Tg]}$), and we have
   \begin{eqnarray}\label{keral}
     \ker(h) & = & \{f\in k[\Tg] \, : \,
     \exists n \, \mm^n.f\subset \ker(\alpha) \} \\ \nonumber
     & = & \{f\in k[\Tg] \, : \,
     \exists n \, f.A_{\geq n}[\Tg]\subset \ker(\alpha) \}.  \nonumber
    \end{eqnarray}
Now it is easy to check that (\ref{ann}) and (\ref{keral})
are the same.
\end{proof}

In what follows we will always suppose that the map $\lambda$ is
generically finite, this implies that the ideal $\ker(h)$ of $k[\Tg]$
is principal. Indeed, it is clear that $\ker(h)$ is a prime ideal of
codimension one since $\lambda$ is supposed to be generically finite. 
As $k[\Tg]$ is a factorial domain, $\ker(h)$ is principal by theorem
1, VII 3.1 \cite{BAC85}.

   
\subsection{Implicitization of a hypersurface without base points}
Here we suppose that both ideals $I$ and $\mm$ of $A$ have the same
radical. This condition is equivalent to say that the variety $V(I)$
is empty in $\Proj(A)$, that is the polynomials $f_1,\ldots f_n$
have no common roots in $\Proj(A)$. Under this condition proposition
\ref{acycwithoutbp} shows that the approximation complex $\Zc$ is
acyclic, since $\depth_\mm(A)=n-1$, and hence gives the exact complex
$$0\rightarrow \Zc_{n-1} \xrightarrow{d_\Tg} \ldots
\xrightarrow{d_\Tg} \Zc_1 \xrightarrow{d_\Tg} A[\Tg] \rightarrow
\Sym_A(I) \rightarrow 0.$$
Recall that $\Zc_i$ and $\Sym_A(I)$ are bigraded. We denote
respectively by $(\Zc_i)_\nu$ and $\Sym_A(I)_\nu$ the
graded part of $\Zc_i$ and $\Sym_A(I)$ corresponding to the graduation
of $A$. We have the following theorem:
\begin{theorem}\label{thnobp} Suppose that both ideals $I$ and $\mm$ have the same
  radical. Let  
  $\eta$ be an integer such that $H^0_\mm(\Sym_A(I))_\nu=0$  for 
  all $\nu\geq \eta$, and denote by $H$ the reduced equation (defined up to a
  nonzero multiple in $k$) of the closed image of the map
  $\lambda$. Then the determinant of the degree $\nu$ part of the
  $\Zc$-complex associated to $I$, which is a complex of
  $k[\Tg]$-modules of the form 
$$0\rightarrow (\Zc_{n-1})_{\nu} \rightarrow \ldots \rightarrow
(\Zc_1)_{\nu} \rightarrow A_\nu[\Tg],$$ 
  is exactly $H^{\deg(\lambda)}$, of degree $d^{n-2}$.
\end{theorem}
\begin{proof} First notice that since $I$ and $\mm$ have the same
  radical, $I$ is of linear type outside $V(\mm)$. Indeed for all prime
  $\mathfrak{p}\notin V(\mm)$ we have $I_\mathfrak{p}=A_\mathfrak{p}$
  and hence, by proposition \ref{indepapprox}, the approximation complex
  $\Mc_\mathfrak{p}$ is acyclic. This implies, by proposition \ref{lintype},
  that $I$ is of linear type outside $V(\mm)$. 

  Now let $\nu$ be a fixed integer greater or equal to
  $\eta$. By proposition \ref{annprop} we know that 
  $$\ann_{k[\Tg]}(H_0(\Zc_\bullet)_\nu)=\ker(h)$$
  which is a principal ideal generated by $H$, that we denote hereafter by $\pp$. Moreover
  by proposition \ref{acycwithoutbp} the $\Zc$-complex associated to
  $I$ is acyclic and so, by standard properties of determinants of complexes (see e.g. \cite{KnMu76}), we deduce that
  $$\det((\Zc_\bullet)_\nu)=\div(H_0(\Zc_\bullet)_\nu)=\length((\Sym_A(I)_\nu)_\pp).[H].$$
It remains to prove that $\length((\Sym_A(I)_\nu)_\pp)=\deg(\lambda)$, but this last equality follows from the proof of theorem \ref{degree}. Using the notation of the theorem \ref{degree}, 
as $I$ is assumed to be of linear type outside $V(\mm)$, 
{\small \begin{eqnarray*}
\length((\Sym_A(I)_\nu)_\pp) & = & 
\dim_{(k[\Tg]/\pp)_\pp}\Gamma(\Proj_A(\Sym_A(I)_\pp),\OO_{\Proj_A(\Sym_A(I)_\pp)}(\nu))\\
& = & \dim_{(k[\Tg]/\pp)_\pp}\Gamma(\Proj_A(B_\pp),
\OO_{\Proj_A(B_\pp)}(\nu))\\
& = & \dim_{(k[\Tg]/\pp)_\pp}\Gamma(\Proj_A(B_\pp),
\OO_{\Proj_A(B_\pp)})=\deg(\lambda),
\end{eqnarray*}}

\noindent where the last equality follows from \eqref{iso} (and the five following lines). Finally theorem \ref{degree} shows also that $H^{\deg(\lambda)}$ is of degree $d^{n-2}$.
\end{proof}

\begin{remark}\label{rmsurface}
Notice that theorem \ref{degree} implies that $\lambda$ is always generically
finite if $d\geq 1$ and $\lambda$ is a regular map (i.e.~there is no
base points). 
\end{remark}
We deduce immediately this corollary which is standard when dealing with 
determinants of complexes (see \cite{GKZ94}, appendix A).
\begin{corollary} Under the hypothesis of theorem \ref{thnobp}, 
$H^{\deg(\lambda)}$ is obtained as the gcd of the maximal minors of
the surjective 
$k[\Tg]$-module morphism $$(\Zc_1)_{\nu} \xrightarrow{d_\Tg}
A_\nu[\Tg],$$ for all $\nu \geq \eta$.
\end{corollary}

It is also possible to give an
explicit bound for the integer $\eta$, bound which is here the best
possible (we will see that it is obtained for particular examples).

\begin{proposition}\label{reg1} ($n\geq 3$) Suppose that both ideals $I$ and $\mm$
  have the same 
  radical,   then $$H^0_\mm(\Sym_A(I))_\nu=0 \ \ \forall  \nu\geq
  (n-2)(d-1).$$
\end{proposition}
\begin{proof} From the standard spectral sequences associated to the
  bicomplex
  $$\begin{array}{ccccccccccc}
   0 & \rightarrow & \CC^0_\mm(\Zc_{n}) &  \xrightarrow{d_\Tg} & \ldots &
   \xrightarrow{d_\Tg} & \CC^0_\mm(\Zc_1) & \xrightarrow{d_\Tg} & \CC^0_\mm(\Zc_0) &
   \rightarrow & 0 \\  
    & & \downarrow &  &   &  & \downarrow &  & \downarrow &  &  \\
  0 & \rightarrow & \CC^1_\mm (\Zc_n) &  \xrightarrow{} & \ldots
  & \xrightarrow{} & 
  \CC^1_\mm (\Zc_1) & \xrightarrow{} & \CC^1_\mm (\Zc_0) & \rightarrow & 0 \\
    & & \downarrow &  &  &  & \downarrow &  & \downarrow &
    &  \\
    & & \vdots &  & \vdots &  & \vdots &  & \vdots &
    &  \\
     & & \downarrow &  &  &  & \downarrow &  & \downarrow &
    &  \\
    0 & \rightarrow & \CC^{n-1}_\mm (\Zc_n) &  \xrightarrow{} & \ldots & \xrightarrow{} &
  \CC^{n-1}_\mm (\Zc_1) & \xrightarrow{} & \CC^{n-1}_\mm (\Zc_0) &  \rightarrow & 0 ,\\ 
\end{array}$$
  where the first row is the $\Zc$-complex associated to $I$, and
  columns are the usual \v{C}ech complexes (as in the proof of lemma
  \ref{tech}),  
  it is easy to see that $H^0_\mm(\Sym_A(I))_\nu=0$ for all $\nu$ such
  that $H^i_\mm(\Zc_i)_\nu=0$ for all $i\geq 1$. We consider the
  Koszul complex $K_\bullet(\fg;A)$ associated to the sequence
  $f_1,\ldots,f_n$ over $A$ and denote by $Z_i$ (resp.~$B_i$) its 
  $i$-cycles (resp.~its $i$-boundaries). As
  $\depth_\mm(A)=\depth_I(A)=n-1$, we know that $B_i=Z_i$ for all
  $i \geq 2$. We deduce first that $Z_n$ equals 0 (and hence $\Zc_n=0$),
  and also that $B_{n-1}\simeq A(-d)$. It follows $H^i_\mm(B_{n-1})=0$
  for $i\neq n-1$ and
  $H^{n-1}_\mm(B_{n-1})_\nu=H^{n-1}_\mm(Z_{n-1})_\nu=0$ for all
  $\nu\geq d-n+2$ (recall that $H^{n-1}(A)_\nu=0$ for all $\nu >
  -(n-1)$). Now by iterations from the exact sequences
  $$0\rightarrow B_{i+1}(-d) \rightarrow K_{i+1}(-d) \rightarrow B_{i}
  \rightarrow 0,$$ for $i\geq 2$, we deduce that
  $$H^i_\mm(B_i)_\nu=H^i_\mm(Z_i)_\nu=0, \ \  \forall \ i \geq 2 \ \mathrm{and} \ \forall \ \nu\geq (n-2)d-n+2.$$
  Finally the exact sequence $0\rightarrow Z_1(-d) \rightarrow A(-d)^n
  \rightarrow I 
  \rightarrow 0$ shows that $H^1_\mm(Z_1)_\nu=0$ for all
  $\nu \geq d-n+2$, and the proposition is proved.
\end{proof}

Before examining the case where base points exist, we
look at the two particular situations of curves and surfaces, and link our
result to the so-called method of \emph{moving surfaces} (see
\ref{overcurve} and \ref{oversurface}). 

\subsubsection{Implicitization of curves in $\PP^2_k$ without base points}\label{curvewbp} 
We consider here the particular case $n=3$. We have 
$\mm=(X_1,X_2) \subset A=k[X_1,X_2]$. We suppose that 
the three homogeneous
polynomials $f_1,f_2,f_3$ of the same degree $d\geq
1$ have no base points, that is they have no common factors (in fact base
points in this case are easily ``erased''  by dividing each $f_i$ by the 
gcd of $f_1,f_2,f_3$; however we will see in \ref{curvebp} that
this computation is not necessary). The $f_i$'s define hence a regular map
$$\PP^1_k \xrightarrow{\lambda} \PP^2_k : (X_1:X_2) \mapsto
(f_1:f_2:f_3)(X_1,X_2),$$
whose image is a curve of reduced implicit equation $C$. All the hypothesis
made in the preceding paragraph are valid and it follows that  the
determinant of each complex $(\Zc_\bullet)_\nu$, for all 
$\nu\geq d-1$, is exactly $C^{\deg(\lambda)}$. These complexes are of
the form
$$0\rightarrow (\Zc_2)_\nu \xrightarrow{d_{\Tg}} (\Zc_1)_\nu
\xrightarrow{d_{\Tg}} A_\nu[\Tg].$$
But $\Zc_2\simeq A(-d)[\Tg]$ since $\depth_I(A)\geq 2$, and hence
$(\Zc_2)_\nu=0$ for all $\nu 
\leq d-1$. We deduce that the determinant of the complex
$(\Zc_\bullet)_{d-1}$, which is $C^{\deg(\lambda)}$, is in fact the
determinant of the matrix $(\Zc_1)_{d-1} \xrightarrow{d_{\Tg}}
A_{d-1}[\Tg]$ of $k[\Tg]$-modules, i.e.~it is obtained as the determinant of the first syzygies of $f_1,f_2,f_3$ in degree
$d-1$. This result is 
exactly the method of \emph{moving lines} which we have recalled in
\ref{overcurve}. 

Consider the simple example $f_1=X_1^2$, $f_2=X_1X_2$ and
$f_3=X_2^2$. Applying the method, the matrix of syzygies of degree
$d-1$ is:
$$\begin{pmatrix}
  -T_2 & T_3 \\
  T_1 & T_2
\end{pmatrix},$$
and the implicit equation is hence $T_2^2-T_1T_3$. The method fails
if we try $\nu=d-2$.

\subsubsection{Implicitization of surfaces in $\PP^3_k$ without base
  points}\label{prsurf} We consider 
here the particular case $n=4$. We suppose that the four
polynomials $f_1,f_2,f_3,f_4$ of the same degree $d$ have no common
roots in $\Proj(A)$, where $A$ is here the polynomial ring
$A=k[X_1,X_2,X_3]$. These polynomials define a regular map
$$\PP^2_k \xrightarrow{\lambda} \PP^3_k : (X_1:X_2:X_3) \mapsto
(f_1:f_2:f_3:f_4)(X_1,X_2,X_3),$$
whose image is a surface of reduced implicit equation $S$. Applying the
preceding results we obtain that the determinant of each complex
$(\Zc_\bullet)_\nu$, for all  
$\nu\geq 2(d-1)$, is exactly $S^{\deg(\lambda)}$. These complexes are of
the form
$$0\rightarrow (\Zc_3)_\nu \xrightarrow{d_{\Tg}} (\Zc_2)_\nu
\xrightarrow{d_{\Tg}} (\Zc_1)_\nu 
\xrightarrow{d_{\Tg}} A_\nu[\Tg].$$
As in the case of curves we have also $\Zc_3\simeq A(-d)[\Tg]$, but
here $2(d-1) \geq d$ since $d\geq 1$. Consequently we obtain here the equation $S^{\deg(\lambda)}$ as the product of two
determinants divided by another one (see \cite{GKZ94},
appendix A, to compute the determinant of a complex). We illustrate it
with the two following examples. 

First consider the example given by $f_1=X_1^2$, $f_2=X_2^2$,
$f_3=X_3^2$ and $f_4=X_1^2+X_2^2+X_3^2$. Clearly the implicit equation is $T_1+T_2+T_3-T_4=0$. Applying our method in degree
  $\nu=2(2-1)=2$ we
  compute the three matrices of the complex which are respectively (from
  the right to the left) of size $6\times 9$, $9\times 4$ and $4\times
  1$. We obtain $(T_1+T_2+T_3-T_4)^4$ as the product of two
  determinants of size $6\times 6$ and $1\times 1$ divided by another
  one of size $3\times 3$. If we try the method in degree $\nu=1$ the
  first matrix (on the right) is square of determinant
  $(T_1+T_2+T_3-T_4)^3$ and the others are zero, but $3$ is different from $4$ which is the degree of the parameterization map.

 As another example we consider $f_1=X_1^2X_2$, $f_2=X_2^2X_3$,
 $f_3=X_1X_3^2$ and $f_4=X_1^3+X_2^3+X_3^3$. Applying our method in
 degree $4=2(3-1)$ we find the irreducible implicit equation of degree
 9 as the product of two determinants of size $15\times 15$ and
 $3\times 3$ divided by another one of size $9\times 9$. The method
 fails in degree strictly less than 4.

\begin{remark}
As we have seen, our method gives in general the implicit
 equation as the product of two determinants divided by another one.  
However, we point out that the method of moving surfaces of T. Sederberg gives 
here (i.e.~in the case without base points) the implicit  
equation as a determinant of a square matrix, as we have recalled in
\ref{oversurface}. 
\end{remark}

\subsection{Implicitization of a hypersurface with isolated local complete
  intersection base points}
In this paragraph we suppose that the ideal $I=(f_1,\ldots,f_{n})$ of the polynomial ring 
$A=k[X_1,\ldots,X_{n-1}]$ is a local complete intersection in 
$\Proj(A)$ of codimension 
$n-2$ (see definition \ref{lcidef}), that is the ideal $I$ defines
points in $\PP^{n-2}_k$ which are 
locally generated by a regular sequence. 
Under this condition we have
$\depth_I(A)=n-2 < \depth_\mm(A)=n-1$, so proposition
\ref{lciprop} shows that the $\Zc$ approximation
complex associated to $I$ is acyclic. As for the case without base
points we obtain  the exact complex 
$$0\rightarrow \Zc_{n-1} \xrightarrow{d_\Tg} \ldots
\xrightarrow{d_\Tg} \Zc_1 \xrightarrow{d_\Tg} A[\Tg] \rightarrow
\Sym_A(I) \rightarrow 0.$$
Recalling that $\Sym_A(I)_\nu$ denotes the
graded part of $\Sym_A(I)$ corresponding to the graduation of $A$, we have the
following theorem: 
\begin{theorem}\label{withbp} Suppose that the ideal
  $I=(f_1,\ldots,f_{n})$ is a 
  local complete intersection in $\Proj(A)$ of codimension $n-2$
  such that the map $\lambda$ is generically finite. Let  
  $\eta$ be an integer such that $H^0_\mm(\Sym_A(I))_\nu=0$  for 
  all $\nu\geq \eta$, and denote by $H$ the reduced equation (defined up to a
  nonzero multiple in $k$) of the closed image of the map
  $\lambda$. Then the determinant of the degree $\nu$ part of the
  $\Zc$-complex associated to $I$, which is a complex of
  $k[\Tg]$-modules of the form 
$$0\rightarrow (\Zc_{n-1})_{\nu} \rightarrow \ldots \rightarrow
(\Zc_1)_{\nu} \rightarrow A_\nu[\Tg],$$ 
  is exactly $H^{\deg(\lambda)}$, of degree
  $d^{n-2}-\dim_k\Gamma(T,\OO_T)$ where $T=\Proj(A/I)$.
\end{theorem}
\begin{proof} First notice that since $I$ is a local complete
  intersection in $\Proj(A)$ we deduce that $H_i(\Mc)$ are supported on $V(\mm)$ for
  all $i\geq 1$, and hence that $I$ is of linear type outside
  $V(\mm)$.  Let $\nu$ be a fixed integer greater or equal to
  $\eta$. By proposition \ref{annprop} we know that 
  $$\ann_{k[\Tg]}(H_0(\Zc_\bullet)_\nu)=\ker(h)$$
  which is a principal ideal generated by $H$. Moreover
  by proposition \ref{lciprop} the
  $\Zc$-complex associated to 
  $I$ is acyclic and hence, by the same argument as the one given in the proof of theorem \ref{thnobp}, we deduce that
  $$\det((\Zc_\bullet)_\nu)=H^{\deg(\lambda)}.$$
  In the same way, theorem \ref{degree} shows that the polynomial $H^{\deg{(\lambda})}$ is of
  degree $d^{n-2}-e(T,\Proj(A))$, and
  $e(T,\Proj(A))=\dim_k\Gamma(T,\OO_T)$ since $T$ is locally a
  complete intersection.
\end{proof}
\begin{remark}\label{rksurface}
The hypothesis saying that $\lambda$ is generically finite can be
tested using theorem \ref{degree}. In the particular case $n=4$ the
hypothesis saying that $I$ is a local complete intersection of
codimension 2 can be tested using  proposition \ref{coxschenck}.
\end{remark}

We deduce the standard corollary:
\begin{corollary}
Under the
hypothesis of theorem \ref{withbp}, 
$H^{\deg(\lambda)}$ is obtained as the gcd of the maximal minors of the
surjective $k[\Tg]$-module morphism $$(\Zc_1)_{\nu} \xrightarrow{d_\Tg}
A_\nu[\Tg],$$
for all $\nu \geq \eta$.
\end{corollary}

We now show that the bound
given for the integer $\eta$ in proposition \ref{reg1} is also
available here. 

\begin{proposition}\label{reg2} ($n\geq 3$) Suppose that the ideal
  $I=(f_1,\ldots,f_{n})$ is a 
  local complete intersection in $\Proj(A)$ of codimension $n-2$,
  then $$H^0_\mm(\Sym_A(I))_\nu=0 \ \
  \forall  \nu\geq (n-2)(d-1).$$ 
\end{proposition}
\begin{proof} This proof is quite similar to the proof of proposition
  \ref{reg1}. The same first argument shows 
  that $H^0_\mm(\Sym_A(I))_\nu=0$ for all $\nu$ such
  that $H^i_\mm(\Zc_i)_\nu=0$ for all $i\geq 1$. We consider the
  Koszul complex $K_\bullet(\fg;A)$ associated to the sequence
  $f_1,\ldots,f_n$ over $A$, and denote by $Z_i$ (resp.~$B_i$) its 
  $i$-cycles (resp.~its $i$-boundaries). As
  $\depth_\mm(A)=n-1$ and $\depth_I(A)=n-2$, we know that $B_i=Z_i$ for all
  $i > 2$. We deduce first that $Z_n$ equals 0 (and hence $\Zc_n=0$),
  and second that $B_{n-1}\simeq A(-d)$. It follows $H^i_\mm(B_{n-1})=0$
  for $i\neq n-1$, and
  $H^{n-1}_\mm(B_{n-1})_\nu=H^{n-1}_\mm(Z_{n-1})_\nu=0$ for all
  $\nu\geq d-n+2$ (recall that $H_\mm^{n-1}(A)_\nu=0$ for all $\nu >
  -(n-1)$). Now by iterations from the exact sequences
  $$0\rightarrow B_{i+1}(-d) \rightarrow K_{i+1}(-d) \rightarrow B_{i}
  \rightarrow 0,$$ for $i\geq 2$, we deduce that
  $$H^i_\mm(B_i)_\nu=H^i_\mm(Z_i)_\nu=0 \ \forall \ i \geq 3 \
  \mathrm{and} \ \forall \ \nu\geq (n-3)d-n+2,$$
  and that $H^2_\mm(B_2)_\nu=0$ for all $\nu \geq (n-2)d-n+2$. 
  Also the exact sequence $0\rightarrow Z_1(-d) \rightarrow A(-d)^n
  \rightarrow I 
  \rightarrow 0$ shows that $H^1_\mm(Z_1)_\nu=0$ for all
  $\nu \geq d-n+2$. Finally $H^2_\mm(H_2(\ag;A))=0$ since $H_2(\ag;A)$
  is supported on $V(I)$ which of dimension 1, and the exact sequence
  $0 \rightarrow B_2 \rightarrow Z_2 \rightarrow H_2 \rightarrow 0$
  shows that $H^2_\mm(Z_2)_\nu=0$ for all $\nu$ such that
  $H^2_\mm(B_2)_\nu=0$, that is for all $\nu \geq (n-2)d-n+2$.
\end{proof}

We now focus on the particular cases of interest of curves and
surfaces.

\subsubsection{Implicitization of curves in $\PP^2_k$ with base
  points}\label{curvebp} 
We consider here the particular case $n=3$. We have 
$\mm=(X_1,X_2) \subset A=k[X_1,X_2]$. We suppose that 
the three homogeneous
polynomials $f_1,f_2,f_3$ of the same degree $d\geq
1$ have base points, that is they have a common factor.
All the hypothesis
made in the preceding paragraph are valid and it follows that  the
determinant of each complex $(\Zc_\bullet)_\nu$, for all 
$\nu\geq d-1$, is exactly $C^{\deg(\lambda)}$, where $C$ denotes the
reduced implicit curve. These complexes are of
the form
$$0\rightarrow (\Zc_2)_\nu \xrightarrow{d_{\Tg}} (\Zc_1)_\nu
\xrightarrow{d_{\Tg}} A_\nu[\Tg].$$
Comparing to \ref{curvewbp}, we do not have $\Zc_2\simeq A(-d)[\Tg]$ 
because here $\depth_I(A)=1$, whereas $\depth_I(A)=2$ in
\ref{curvewbp}.  The determinant of the complex
$(\Zc_\bullet)_{d-1}$, which is $C^{\deg(\lambda)}$, is hence generally
obtained as the quotient of two determinants. We can illustrate this
with the simple example of \ref{curvewbp} where we multiply each
equation by $X_1$: $f_1=X_1^3$, $f_2=X_1^2X_2$ and
$f_3=X_1X_2^2$. The first matrix on the right is given by:
$$\begin{pmatrix}
 -T_2 & -T_3 & -T_3 &  0\\
 T_1 &  0 &  T_2 &  -T_3 \\
 0 & T_1 &  0 &  T_2
\end{pmatrix},$$
and the second one is given by:
$$\begin{pmatrix}
 -T_3\\
 T_2\\
 0\\
-T_1
\end{pmatrix}.$$
It follows that the implicit equation is obtained as the quotient
$$\frac{\left|\begin{array}{ccc}
 -T_2 & -T_3 & -T_3 \\
 T_1 &  0 &  T_2  \\
 0 & T_1 &  0
 \end{array}\right|}{|-T_1|}=-T_2^2+T_1T_3.$$

\subsubsection{Implicitization of surfaces in $\PP^3_k$ with lci base points}
We suppose here that $n=4$. As pointed out in remark \ref{rksurface}, 
it is here possible to know the degree of the implicit equation we are
looking for and to test if the base points are locally complete
intersection or not. Under the hypothesis of theorem
\ref{withbp} we obtain here again 
that the implicit equation is the determinant of the
complex $(\Zc_\bullet)_\nu$, for all $\nu\geq 2(d-1)$. As in the case where
there is no base points, these complexes are of the form   
$$0\rightarrow (\Zc_3)_\nu \xrightarrow{d_{\Tg}} (\Zc_2)_\nu
\xrightarrow{d_{\Tg}} (\Zc_1)_\nu 
\xrightarrow{d_{\Tg}} A_\nu[\Tg].$$
with $\Zc_3\simeq A(-d)[\Tg]$. Consequently we deduce the
implicit equation as the product of two
determinants divided by another one. 

Notice that recently it has been proved in \cite{BCD02} that the
method of moving 
surfaces works with local complete intersection base points which
satisfy additional conditions (see \ref{oversurface}).
 The algorithm we propose here avoids these technical 
hypothesis. To illustrate it we end 
this section with the following example taken from \cite{BCD02}, for 
which the (improved) moving surfaces method failed. We take $f_1=X_1X_3^2$,
$f_2=X_2^2(X_1+X_3)$, $f_3=X_1X_2(X_1+X_3)$ and 
  $f_4=X_2X_3(X_1+X_3)$. The ideal $I=(f_1,f_2,f_3,f_4)$, which is not 
  saturated, is a local complete
  intersection in $\PP^2$ of codimension 2  defining 6 points (counted with
  multiplicity). Applying our method  in degree $2(3-1)=4$, we obtain the implicit equation as the quotient
  $\frac{\Delta_0\Delta_2}{\Delta_1}$, where $\Delta_0$ is the
  determinant of the matrix 
  $${\tiny \left(\begin{array}{ccccccccccccccc}
0 & 0 & 0 & 0 & 0 & 0 & T_2 & 0 & 0 & 0 & 0 & 0 & 0 & 0 & 0\\
T_1 &0 &0 &0 &0 &0 &-T_3 &T_2 &-T_4 &0 &0 &0 &0 &0 &0\\
-T_4 &0 &0 &0 &0 &0 &0 &0 &T_2 &0 &0 &0 &0 &0 &0\\
0 &T_1 &0 &0 &0 &0 &0 &-T_3 &0 &T_2 &-T_4 &0 &0 &0 &0\\
0 &-T_4 &T_1 &0 &0 &0 &0 &0 &0 &0 &T_2 &-T_4 &0 &0 &0\\
T_4 &0 &-T_4 &0 &0 &0 &0 &0 &0 &0 &0 &T_2 &0 &0 &0\\
0 &0 &0 &T_1 &0 &0 &0 &0 &0 &-T_3 &0 &0 &T_2 &0 &0\\
0 &0 &0 &-T_4 &T_1 &0 &0 &0 &0 &0 &0 &0 &0 &0 &0\\
0 &T_4 &0 &0 &-T_4 &T_1 &0 &0 &0 &0 &0 &0 &0 &-T_4 &0\\
-T_4 &0 &T_4 &0 &0 &-T_4 &0 &0 &0 &0 &0 &0 &0 &T_2 &0\\
0 &0 &0 &0 &0 &0 &0 &0 &0 &0 &0 &0 &-T_3 &0 &0\\
0 &0 &0 &T_1 &0 &0 &0 &0 &0 &0 &0 &0 &0 &0 &0\\
0 &-T_1 &0 &0 &T_1 &0 &0 &0 &0 &0 &0 &0 &0 &0 &0\\
T_1 &0 &-T_1 &0 &0 &T_1 &0 &0 &0 &0 &0 &0 &0 &0 &-T_4\\
0 &0 &0 &0 &0 &0 &0 &0 &0 &0 &0 &0 &0 &0 &T_2
\end{array}\right),}$$
$\Delta_1$ is the determinant of 
$${\tiny \left(\begin{array}{ccccccccccccccc}
 0 &0 &0 &T_1 &0 &0 &0 &0 &0 &0 &0 &-T_3 &0 &0 &0\\
 0 &T_4 &0 &-T_4 &T_1 &0 &0 &0 &0 &0 &0 &0 &0 &0 &0\\
 0 &0 &0 &0 &0 &0 &0 &0 &0 &0 &0 &0 &0 &-T_3 &0\\
 0 &-T_1 &0 &T_1 &0 &0 &0 &0 &0 &0 &0 &0 &0 &0 &0\\
 T_1 &0 &-T_1 &0 &T_1 &0 &0 &0 &0 &0 &0 &0 &0 &0 &0\\
 0 &0 &0 &0 &0 &0 &0 &0 &T_2 &0 &0 &0 &0 &0 &0\\
 T_4 &0 &0 &0 &0 &T_1 &0 &0 &0 &T_2 &-T_4 &0 &0 &0 &0\\
 0 &0 &0 &0 &0 &-T_4 &0 &0 &0 &0 &T_2 &0 &0 &0 &0\\
 0 &T_4 &0 &0 &0 &0 &T_1 &0 &0 &0 &0 &T_2 &0 &0 &0\\
-T_4 &0 &0 &0 &0 &0 &-T_4 &T_1 &0 &0 &0 &0 &-T_4 &0 &0\\
 0 &0 &0 &0 &0 &T_4 &0 &-T_4 &0 &0 &0 &0 &T_2 &0 &0\\
 0 &-T_1 &0 &0 &0 &0 &0 &0 &0 &0 &0 &0 &0 &T_2 &0\\
 T_1 &0 &0 &0 &0 &0 &T_1 &0 &0 &0 &0 &0 &0 &0 &0\\
 0 &0 &0 &0 &0 &-T_1 &0 &T_1 &0 &0 &0 &0 &0 &0 &-T_4\\
 0 &0 &0 &0 &0 &0 &0 &0 &0 &0 &0 &0 &0 &0 &T_2
\end{array}\right),}$$
and $\Delta_3$ is the determinant of
$${\tiny \left(\begin{array}{ccc}
 0 &-T_4 & T_1\\
 0 & T_1 & 0\\
-T_1 & 0 & T_1
\end{array}\right).}$$
We deduce that the desired equation is
$T_1T_2T_3+T_1T_2T_4-T_3T_4^2$. Point out that, trying our method 
 empirically, we find the latter equation
by applying theorem \ref{withbp} with $\nu=3$, and even $\nu=2$. For
the case $\nu=3$ the equation is obtained as a similar quotient
$\frac{\Delta_0\Delta_2}{\Delta_1}$, where here $\Delta_0$ is of size
$10\times 10$, $\Delta_1$ of size $8\times 8$, and $\Delta_2$ of size
$1\times 1$. For the case $\nu=2$, the equation is obtained as a
quotient $\frac{\Delta_0}{\Delta_1}$ (the third map of the
$\Zc$-complex in this degree degenerates to zero), where $\Delta_0$ is 
of size $6\times 6$ and $\Delta_1$ of size $3\times 3$. 

\appendix
\section{Kravitsky's formula and anisotropic resultant}\label{app}

Let $k$ be a commutative ring, 
$X$ and $Y$ be two indeterminates, and $P,Q,R \in k[X,Y]_d$ be three
homogeneous polynomials of same given degree $d\geq 1$. Introducing
three new indeterminates $u,v,w$, we may consider the following
$d\times d$ matrix with coefficients in $k[u,v,w]$:
$$\Omega:=u\Bez(Q,R)+v\Bez(R,P)+w\Bez(P,Q).$$

Let $Z$ be another new indeterminate. Setting the weight of $X,Y,Z$
respectively to $1,1,d$, the polynomials $P-uZ,Q-vZ,R-wZ$ are isobaric 
of weight $d$; we can hence consider the anisotropic
resultant (see \cite{Jou91} and \cite{Jou96})

$$\mathcal{R}:= {}^a \Res(P-uZ,Q-vZ,R-wZ) \in
k[u,v,w].$$

\noindent \textbf{Proposition} {\it With the preceding notations, we have the
equality
$$\mathcal{R}=(-1)^{\frac{d(d+1)}{2}}\det(\Omega)$$
in $k[u,v,w]$. }

\begin{proof}
By specialization, it is sufficient to prove this proposition in case
$P,Q,R$ have indeterminate coefficients and
$k=\ZZ[\mathrm{coeff}(P,Q,R)]$, that we suppose hereafter. The
generic polynomials $P-uZ,Q-vZ,R-wZ$ are isobar of weight $d$, and $d$
 is a multiple of the lcm of the weights of $X,Y,Z$. In this way
$\mathcal{R}$ is \emph{irreducible}, and normalized by the condition
$${}^a \Res(X^d,Y^d,Z)=1.$$
We denote $A:=\ZZ[\mathrm{coeff}(P,Q,R)][u,v,w]$ and
$B:=A[X,Y,Z]/(P-uZ,Q-vZ,R-wZ)$. We know that the ideal $H^0_\mm(B)$,
where $\mm=(X,Y,Z)$, is
prime, and that, by definition, $\mathcal{R}$ is a generator of the
ideal $H^0_\mm(B)_0$ of $A$. For all $L\in k[X,Y]$ we denote
$\tilde{L}:=L(X,1) \in k[X]$, and it follows that the ring
$$\tilde{B}:=A[X,Z]/(\tilde{P}-uZ,\tilde{Q}-vZ,\tilde{R}-wZ)$$ has no
zero divisors, and that
$$H^0_\mm(B)_0=A\cap (\tilde{P}-uZ,\tilde{Q}-vZ,\tilde{R}-wZ) \
\mathrm{in} \ A[X,Z].$$
Finally, we know that
$$\Res(P-uZ^d,Q-vZ^d,R-wZ^d)= \mathcal{R}^d,$$
where $\Res$ denotes the classical resultant of three homogeneous
polynomials, which shows that $\mathcal{R}$ is homogeneous of degree
$3d^2/d=3d$ in the coefficients of the polynomials
$P-uZ,Q-vZ,R-wZ$.

Now introduce new indeterminates $T_1,\ldots,T_d$, and denote by
$$\Delta(T_1,\ldots,T_d):=(T^{d-j})_{1\leq i,j \leq d}$$
the Vandermonde's matrix. We have clearly
\begin{equation}\label{vander}
\Delta(T_1,\ldots,T_d)\circ \Omega \circ
\left(\begin{array}{c}
X^{d-1}\\
X^{d-2}\\
\vdots\\
1
\end{array}\right)=
\left(\begin{array}{c}
a_1\\
a_2\\
\vdots\\
a_d
\end{array}\right),
\end{equation}
with, for all $1\leq i \leq d$,
\begin{eqnarray*}
  (T_i-X)a_i & = &
  u(\tilde{Q}(T_i)\tilde{R}(X)-\tilde{R}(T_i)\tilde{Q}(X))\\
  & + & v(\tilde{R}(T_i)\tilde{P}(X)-\tilde{P}(T_i)\tilde{R}(X))\\
  & + & w(\tilde{P}(T_i)\tilde{Q}(X)-\tilde{Q}(T_i)\tilde{P}(X)),
\end{eqnarray*}
which can be rewritten
\begin{eqnarray*}
  (T_i-X)a_i & = &
  u(\tilde{Q}(T_i)(\tilde{R}(X)-wZ)-\tilde{R}(T_i)(\tilde{Q}(X)-vZ))\\
  & + & v(\tilde{R}(T_i)(\tilde{P}(X)-uZ)-\tilde{P}(T_i)(\tilde{R}(X)-wZ))\\
  & + & w(\tilde{P}(T_i)(\tilde{Q}(X)-vZ)-\tilde{Q}(T_i)(\tilde{P}(X)-uZ)).
\end{eqnarray*}
This shows that $(T_i-X)a_i$ is in the ideal of 
$A[T_1,T_2,\ldots,T_d,X,Z]$ generated by the polynomials
$\tilde{P}-uZ,\tilde{Q}-vZ,\tilde{R}-wZ$. Composing \eqref{vander} on
the left by the adjoint matrix of $\Delta\circ \Omega$, we deduce that,
for all $0\leq l \leq d-1$,
$$(\prod_{1\leq i \leq d}(T_i-X)\prod_{1\leq i<j\leq d}
  (T_i-T_j))\det(\Omega)X^l \in
(\tilde{P}-uZ,\tilde{Q}-vZ,\tilde{R}-wZ).$$
As $(\prod_{1\leq i \leq d}(T_i-X)\prod_{1\leq i<j\leq d}
  (T_i-T_j))$ is not a zero divisor in $B[T_1,\ldots,T_d]$ (by the
Dedekind-Mertens lemma), it follows in particular (case $l=0$) that
$$\det(\Omega) \in (\tilde{P}-uZ,\tilde{Q}-vZ,\tilde{R}-wZ),$$
and thus $\mathcal{R}$ divides $\det(\Omega)$ in the ring $A$. The
entries of the matrix $\Omega$ are homogeneous of degree 3 with
respect to $u,v,w$ and $\mathrm{coeff}(P,Q,R)$, thus 
$\det(\Omega)$ is homogeneous of degree $3d$ with respect to those
coefficients.  Consequently, there exists $c\in \ZZ$ such that
$\det(\Omega)=c\mathcal{R}$. We can determine $c$ with the
specialization $P\mapsto X^d$,  $Q\mapsto Y^d$,   $R\mapsto 0$,
$u\mapsto 0$, $v\mapsto 0$ and $w\mapsto w$. On one hand the matrix $\Omega$
specializes in $w\Bez(X^d,Y^d)$, and hence $\det(\Omega)$ specializes
in
$$(-1)^{\frac{d(d-1)}{2}}w^d\Res(X^d,Y^d)=(-1)^{\frac{d(d-1)}{2}}w^d,$$
and on the other hand $\mathcal{R}$ specializes in
$${}^a\Res(X^d,Y^d,-wZ)=(-1)^dw^d.{}^a\Res(X^d,Y^d,Z)=(-1)^dw^d,$$
following the properties of the anisotropic resultant. The comparison
of these two last equalities gives the proposition.
\end{proof}

\section*{Acknowledgement} The authors are grateful to the referee for useful comments and remarks.


\begin{thebibliography}{10}
\expandafter\ifx\csname url\endcsname\relax
  \def\url#1{\texttt{#1}}\fi
\expandafter\ifx\csname urlprefix\endcsname\relax\def\urlprefix{URL }\fi

\bibitem{ArSe01}
F.~Ari\`es, R.~Senoussi, An implicitization algorithm for rational surfaces
  with no base points, J. of Symbolic Computation 31 (2001) 357--365.

\bibitem{BAC85}
N.~Bourbaki, Alg\`ebre Commutative, Masson-Dunod, 1985.

\bibitem{BrHe93}
W.~Bruns, J.~Herzog, Cohen-{M}acaulay rings, Cambridge studies in advanced
  mathematics 39.

\bibitem{BuEi73}
D.~A. Buchsbaum, D.~Eisenbud, What makes a complex exact?, Journal of Algebra
  25 (1973) 259--268.

\bibitem{BusPhD}
L.~Bus\'e, \'Etude du r\'esultant sur une vari\'et\'e alg\'ebrique, PhD thesis,
  University of Nice, 2001.

\bibitem{Bus01}
L.~Bus\'e, Residual resultant over the projective plane and the implicitization
  problem, Proceedings ISSAC2001  (2001) 48--55.

\bibitem{BCD02}
L.~Bus\'e, D.~Cox, C.~D'Andrea, Implicitization of surfaces in $\mathbb{P}^3$
  in the presence of base points, to appear in J. of Algebra and its 
  Applications. 

\bibitem{Cox01}
D.~A. Cox, Equations of parametric curves and surfaces via syzygies,
  Contemporary Mathematics 286 (2001) 1--20.

\bibitem{CoSc01}
D.~A. Cox, H.~Schenck, Local complete intersections in $\mathbb{P}^2$ and
  Koszul syzygies, Preprint mathAG/0110097.

\bibitem{Dan01}
C.~D'{A}ndrea, Resultants and moving surfaces, J. of Symbolic Computation 31 (2001) 585--602.

\bibitem{Eis94}
D.~Eisenbud, {C}ommutative {A}lgebra with a view toward {A}lgebraic {G}eometry,
  Vol. 150 of Graduate Texts in Math., Springer-Verlag, 1994.

\bibitem{GKZ94}
I.~Gelfand, M.~Kapranov, A.~Zelevinsky, Discriminants, {R}esultants and
  {M}ultidimensional {D}eterminants, Birkh{\"{a}}user, Boston-Basel-Berlin,
  1994.

\bibitem{GrDi71}
A.~Grothendieck, J.~Dieudonn\'e, \'El\'ements de g\'eom\'etrie alg\'ebrique I,
  Springer-Verlag Berlin Heidelberg New York, 1978.

\bibitem{Har77}
R.~Hartshorne, Algebraic Geometry, Springer-Verlag, 1977.

\bibitem{HSV82}
J.~Herzog, A.~Simis, W.~Vasconcelos, Approximation complexes of blowing-up
  rings, J. of Algebra 74 (1982) 466--493.

\bibitem{HSV83}
J.~Herzog, A.~Simis, W.~Vasconcelos, Approximation complexes of blowing-up
  rings II, J. of Algebra 82 (1983) 53--83.

\bibitem{Jou91}
J.-P. Jouanolou, Le formalisme du r\'esultant, Adv. in Math. 90~(2) (1991)
  117--263.

\bibitem{Jou96}
J.-P. Jouanolou, R\'esultant anisotrope: Compl\'ements et applications, The
  electronic journal of combinatorics 3~(2).

\bibitem{Jou97}
J.-P. Jouanolou, Formes d'inertie et r\'esultant: un formulaire, Adv. in Math.
  126~(2) (1997) 119--250.

\bibitem{KnMu76}
F.~Knudsen, D.~Mumford, The projectivity of the moduli space of stable curves.
  I: Preliminaries on {D}et and {D}iv, Math. Scand. 39 (1976) 19--55.

\bibitem{LKMV95}
M.~S. Livsic, N.~Kravitsky, A.~S. Markus, V.~Vinnikov, Theory of commuting
  nonselfadjoint operators, Mathematics and its Applications, 332. Kluwer
  Academic Publishers Group, Dordrecht, 1995.

\bibitem{MuFo82}
D.~Mumford, J.~Fogarty, Geometric invariant theory - Second edition,
  Springer-Verlag, 1982.

\bibitem{SeCh95}
T.~Sederberg, F.~Chen, Implicitization using moving curves and surfaces,
  Proceedings of SIGGRAPH  (1995) 301--308.

\bibitem{SiVa81}
A.~Simis, W.~Vasconcelos, The syzygies of the conormal module, American J.
  Math. 103 (1981) 203--224.

\bibitem{Vas94}
W.~Vasconcelos, Arithmetic of Blowup Algebras, Vol. 195 of London Mathematical
  Society Lecture Note Series, Cambridge University Press, 1994.

\end{thebibliography}
\end{document}